# DATA-DRIVEN RATE-OPTIMAL SPECIFICATION TESTING IN REGRESSION MODELS[1]

By Emmanuel Guerre and Pascal Lavergne

*LSTA Paris 6 and University of Toulouse GREMAQ and INRA*

We propose new data-driven smooth tests for a parametric regression function. The smoothing parameter is selected through a new criterion that favors a large smoothing parameter under the null hypothesis. The resulting test is adaptive rate-optimal and consistent against Pitman local alternatives approaching the parametric model at a rate arbitrarily close to $1/\sqrt{n}$. Asymptotic critical values come from the standard normal distribution and the bootstrap can be used in small samples. A general formalization allows one to consider a large class of linear smoothing methods, which can be tailored for detection of additive alternatives.

**1. Introduction.** Consider $n$ observations $(Y_i, X_i)$ in $\mathbb{R} \times \mathbb{R}^p$ and the heteroscedastic regression model with unknown mean $m(\cdot)$ and variance $\sigma^2(\cdot)$,

$$Y_i = m(X_i) + \varepsilon_i, \qquad \mathbb{E}[\varepsilon_i | X_i] = 0 \quad \text{and} \quad \text{Var}[\varepsilon_i | X_i] = \sigma^2(X_i).$$

We want to test the hypothesis that the regression belongs to some parametric family $\{\mu(\cdot; \theta); \theta \in \Theta\}$, that is,

$$(1.1) \qquad H_0: m(\cdot) = \mu(\cdot; \theta) \qquad \text{for some } \theta \in \Theta.$$

Tests of $H_0$ are called lack-of-fit tests or specification tests. Based on smoothing techniques, many consistent tests of $H_0$ have been proposed, the so-called smooth tests; see Hart (1997) for a review. A fundamental issue is the choice of the smoothing parameter. Since this is a model selection problem, Eubank and Hart (1992), Ledwina (1994), Hart [(1997), Chapter 7] and Aerts, Claeskens and Hart (1999, 2000), among others, have proposed use of criteria developed by Akaike (1973) and Schwarz (1978). However, these criteria

Received July 2003; revised April 2004.

[1]Supported by LSTA and INRA.

*AMS 2000 subject classifications.* Primary 62G10; secondary 62G08.

*Key words and phrases.* Hypothesis testing, nonparametric adaptive tests, selection methods.

are tailored for estimation but not for testing purposes. Hence, they do not yield adaptive rate-optimal tests, that is, tests that detect alternatives of unknown smoothness approaching the null hypothesis at the fastest possible rate when the sample size grows; see Spokoiny (1996).

Many adaptive rate-optimal specification tests are based on the maximum approach, which consists of choosing as a test statistic the maximum of Studentized statistics associated with a sequence of smoothing parameters. This approach is used for testing the white noise model with normal errors by Fan (1996) and for testing a linear regression model with normal errors by Fan and Huang (2001) and Baraud, Huet and Laurent (2003), who extend the maximum approach. Further work on the linear model includes Spokoiny (2001) under homoscedastic errors and Zhang (2003) under heteroscedastic errors. Finally, Horowitz and Spokoiny (2001) deal with the general case of a nonlinear model with heteroscedastic errors.

We reconsider the model selection approach to propose a new test with some distinctive features. First, our data-driven choice of the smoothing parameter relies on a specific criterion tailored for testing purposes. This yields an adaptive rate-optimal test. Second, the criterion favors a baseline statistic under the null hypothesis. This results in a simple asymptotic distribution for our statistic and in bounded critical values for our test. By contrast, in the maximum approach, critical values diverge and must practically be evaluated by simulation for any sample size. The computational burden of this task can be heavy for a large sample size and a large number of statistics. Moreover, diverging critical values are expected to yield some loss of power compared to our test. In particular, from an asymptotic viewpoint, our test detects local Pitman alternatives converging to the null at a faster rate than the ones detected by a maximum test. In small samples, our simulations show that our test has better power than a maximum test against irregular alternatives.

In our work we allow for a nonlinear parametric regression model with mutidimensional covariates, nonnormal errors and heteroscedasticity of unknown form. In Section 2 we describe the specific aspects of our testing procedure. In Section 3 we detail the practical construction of the test statistic for three types of smoothing procedures. Then we give our assumptions and main results, which concern the null asymptotic behavior of the test, adaptive rate-optimality, and detection of Pitman local alternatives. In Section 4 we prove the validity of a bootstrap method and compare the small sample performances of our test with a maximum test through a simulation experiment. In Section 5 we extend our results to general linear smoothing methods. Finally, we propose a test whose power against additive alternatives is not affected by the curse of dimensionality. Proofs are given in Section 6.

**2. Description of the procedure.** Consider a collection $\{\widehat{T}_h, h \in \mathcal{H}_n\}$ of asymptotically centered statistics which measures the lack-of-fit of the null parametric model. The index $h$ is a smoothing parameter, chosen in a discrete grid whose cardinality grows with the sample size $n$; see our examples in the next section. A maximum test rejects $H_0$ when $\max_{h \in \mathcal{H}_n} \widehat{T}_h/\widehat{v}_h \geq z_\alpha^{\max}$, where $\widehat{v}_h$ estimates the asymptotic null standard deviation of $\widehat{T}_h$. A test in the spirit of Baraud, Huet and Laurent ([2003](#)) rejects the null if $\widehat{T}_h \geq \widehat{v}_h z_\alpha(h)$ for some $h$ in $\mathcal{H}_n$ or, equivalently, if $\max_{h \in \mathcal{H}_n}(\widehat{T}_h/\widehat{v}_h - z_\alpha(h)) > 0$, where the critical values are chosen to get an asymptotic $\alpha$-level test, a difficult issue in practice. Setting $z_\alpha(h) = z_\alpha^{\max}$ yields a maximum test. Because the number $h$ increases with $n$, $z_\alpha^{\max}$ diverges.

On an informal basis, our approach favors a baseline statistic $\widehat{T}_{h_0}$ with lowest variance among the $\widehat{T}_h$. In practice, $\widehat{T}_{h_0}$ can be designed to yield high power against parametric or regular alternatives that are of primary interest for the statistician. However, this statistic may not be powerful enough against nonparametric or irregular alternatives. We then propose to combine this baseline statistic with the other statistics $\widehat{T}_h$ in the following way. Let $\widehat{v}_{h,h_0}$ be some positive estimators of the asymptotic null standard deviation of $\widehat{T}_h - \widehat{T}_{h_0}$. We select $h$ as

$$
\begin{aligned}
(2.1) \qquad \widetilde{h} &= \arg\max_{h \in \mathcal{H}_n}\{\widehat{T}_h - \gamma_n \widehat{v}_{h,h_0}\} \\
&= \arg\max_{h \in \mathcal{H}_n}\{\widehat{T}_h - \widehat{T}_{h_0} - \gamma_n \widehat{v}_{h,h_0}\} \qquad \text{where } \gamma_n > 0.
\end{aligned}
$$

Our test is

$$
(2.2) \qquad \text{Reject } H_0 \text{ when } \widehat{T}_{\widetilde{h}}/\widehat{v}_{h_0} \geq z_\alpha,
$$

where $z_\alpha$ is the quantile of order $(1-\alpha)$ of a standard normal.

The distinctive features of our approach are as follows. First, our criterion penalizes each statistic by a quantity proportional to its standard deviation, while the criteria reviewed in Hart ([1997](#)) use a larger penalty proportional to the variance. Second, the data-driven choice of the smoothing parameter favors $h_0$ under the null hypothesis. Indeed, since $\widehat{T}_h - \widehat{T}_{h_0}$ is of order $\widehat{v}_{h,h_0}$ under $H_0$, $\widetilde{h} = h_0$ asymptotically under $H_0$ if $\gamma_n$ diverges fast enough; see Theorem 1 below. Hence, the null limit distribution of the test statistic is the one of $\widehat{T}_{h_0}/\widehat{v}_{h_0}$, that is, the standard normal, and the resulting test has bounded critical values. Third, our selection procedure allows us to choose the standardization $\widehat{v}_{h_0}$. We could use $\widehat{v}_{\widetilde{h}}$ instead, which also gives an asymptotic $\alpha$-level test since $\widetilde{h} = h_0$ asymptotically under $H_0$. But, because $\widehat{v}_h \geq \widehat{v}_{h_0}$ asymptotically for any admissible $h$, our standardization gives a larger critical region under the alternative. This increases power at no cost from an asymptotic viewpoint; see Fan ([1996](#)) for a similar device in wavelet thresholding tests. Our simulation results show that this effect is already large in small samples. By contrast, the maximum approach systematically downweights the statistic $\widehat{T}_h$ with its standard deviation.

Third, compared to a test using a single statistic, our test inherits the power properties of each of the $\widehat{T}_h$, up to a term $\gamma_n \widehat{v}_{h,h_0}$. Indeed, the definition of $\widetilde{h}$ yields

$$\widehat{T}_{\widetilde{h}} = \max_{h \in \mathcal{H}_n}(\widehat{T}_h - \gamma_n \widehat{v}_{h,h_0}) + \gamma_n \widetilde{v}_{\widetilde{h},h_0} \geq \widehat{T}_h - \gamma_n \widehat{v}_{h,h_0} \qquad \text{for any } h \in \mathcal{H}_n.$$

As a consequence, a lower bound for the power of the test is

$$(2.3) \qquad \mathbb{P}(\widehat{T}_{\widetilde{h}} \geq \widehat{v}_{h_0} z_\alpha) \geq \mathbb{P}(\widehat{T}_h \geq \widehat{v}_{h_0} z_\alpha + \gamma_n \widehat{v}_{h,h_0}) \qquad \text{for any } h \text{ in } \mathcal{H}_n.$$

Using a penalty proportional to a standard deviation yields a better power bound than the selection criteria reviewed in Hart (1997). A suitable choice of the smoothing parameter in the latter power bound allows us to establish the adaptive rate-optimality of the test; see Theorem 2 below and the following discussion. Fourth, combining the $\widehat{T}_h$ with our selection procedure gives a more powerful test than using the baseline statistic $\widehat{T}_{h_0}$. Indeed, since $\widehat{v}_{h_0,h_0} = 0$, a noteworthy implication of (2.3) is

$$(2.4) \qquad \mathbb{P}(\widehat{T}_{\widetilde{h}} \geq \widehat{v}_{h_0} z_\alpha) \geq \mathbb{P}(\widehat{T}_{h_0} \geq \widehat{v}_{h_0} z_\alpha).$$

Theorem 3 below uses the latter inequality to study detection of Pitman local alternatives approaching the null at a faster rate than in Horowitz and Spokoiny (2001).

**3. Main results.** For any integer $q$ and any $x \in \mathbb{R}^q$, $|x| = \max_{1 \leq i \leq q} |x_i|$. For real deterministic sequences, $a_n \asymp b_n$ means that $a_n$ and $b_n$ have the same exact order, that is, there is a $C > 1$ with $1/C \leq a_n/b_n \leq C$ for $n$ large enough. For real random variables, $A_n \asymp_{\mathbb{P}} B_n$ means that $\mathbb{P}(1/C \leq A_n/B_n \leq C)$ goes to 1 when $n$ grows. In such statements, uniformity with respect to a variable means that $C$ can be chosen independently of it. A sequence $\{m_n(\cdot)\}_{n \geq 1}$ is equicontinuous if, for any $\epsilon > 0$, there is an $\eta > 0$ such that $\sup_{n \geq 1} |m_n(x) - m_n(x')| \leq \epsilon$ for all $x$, $x'$ with $|x - x'| \leq \eta$.

3.1. *Construction of the statistics and assumptions.* Let $\widehat{\theta}_n$ be the nonlinear least-squares estimator of $\theta$ in model (1.1), that is,

$$(3.1) \qquad \widehat{\theta}_n = \arg\min_{\theta \in \Theta} \sum_{i=1}^n (Y_i - \mu(X_i; \theta))^2,$$

with an appropriate convention in case of ties. A typical statistic $\widehat{T}_h$ is an estimator of the mean-squared distance of the regression function from the parametric model

$$(3.2) \qquad \min_{\theta \in \Theta} \sum_{i=1}^n (m_n(X_i) - \mu(X_i; \theta))^2.$$

From the estimated parametric residuals $\widehat{U}_i = Y_i - \mu(X_i; \widehat{\theta}_n) = m(X_i) - \mu(X_i; \widehat{\theta}_n) + \varepsilon_i, i = 1, \ldots, n$, we can estimate the departure from the parametric regression using a leave-one-out linear nonparametric estimator $\widehat{\delta}_h(X_i) = \sum_{j=1, j \neq i}^n \nu_{ij}(h)\widehat{U}_j$ based on some weights $\nu_{ij}(h)$ with smoothing parameter $h$. Then (3.2) can be estimated as

$$(3.3) \qquad \widehat{T}_h = \sum_{i=1}^n \widehat{U}_i \widehat{\delta}_h(X_i) = \sum_{1 \leq i \neq j \leq n} \frac{\nu_{ij}(h) + \nu_{ji}(h)}{2} \widehat{U}_i \widehat{U}_j = \widehat{U}' W_h \widehat{U},$$

where $\widehat{U} = [\widehat{U}_1, \ldots, \widehat{U}_n]'$ and the generic element of $W_h$ is $w_{ij}(h) = (\nu_{ij}(h) + \nu_{ji}(h))/2$ for $i \neq j$ and $w_{ii}(h) = 0$. Such a $\widehat{T}_h$ is asymptotically normal under $H_0$; see, for example, de Jong (1987). Examples 1a and 1b come from projection methods, while Example 2 builds on kernel smoothing.

EXAMPLE 1A (Regression on multivariate polynomial functions). Let $\psi_k(x) = \prod_{\ell=1}^p x_\ell^{k_\ell}$, for $k \in \mathbb{N}^p$ with $|k| = \max_{l=1,\ldots,p} k_l \leq 1/h$. Let $\Psi_h = [\psi_k(X_i), |k| \leq 1/h, i = 1, \ldots, n]$ and $P_h = \Psi_h(\Psi_h'\Psi_h)^{-1}\Psi_h'$ be the $n \times n$ orthogonal projection matrix onto the linear subspace of $\mathbb{R}^n$ spanned by $\Psi_h$. The matrix $W_h$ is obtained from $P_h$ by setting its diagonal elements to zero.

EXAMPLE 1B (Regression on piecewise polynomial functions). Under the assumption that the support of $X$ is $[0,1]^p$, we consider piecewise polynomial functions of fixed order $\bar{q}$ over bins $I_k(h) = \prod_{\ell=1}^p [k_\ell h, (k_\ell + 1)h)$, $k = (k_1, \ldots, k_p)$, $k_\ell = 0, \ldots, (1/h) - 1$. These functions write

$$\psi_{qkh}(x) = \prod_{\ell=1}^p x_\ell^{q_\ell} \mathbb{I}(x \in I_k(h)),$$

$$0 \leq |q| = \max_{1 \leq \ell \leq p} q_\ell \leq \bar{q}, 1 \leq |k| = \max_{1 \leq \ell \leq p} k_\ell \leq 1/h.$$

The particular choice $\bar{q} = 0$ corresponds to the regressogram. The matrix $W_h$ is constructed as in Example 1a.

EXAMPLE 2 (Kernel smoothing). Consider a continuous, nonnegative, symmetric and bounded kernel $K(\cdot)$ from $\mathbb{R}^p$ that integrates to 1 and has a positive integrable Fourier transform. These conditions hold for products of the triangular, normal, Laplace or Cauchy kernels. Define $K_h(x) = K(x_1/h, \ldots, x_p/h)$. We consider

$$\widehat{T}_h = \sum_{1 \leq i \neq j \leq n} \frac{1}{(n-1)h^p} \widehat{U}_i \frac{K_h(X_i - X_j)}{\sqrt{\widehat{f}_h(X_i)\widehat{f}_h(X_j)}} \widehat{U}_j$$

with $\widehat{f}_h(X_i) = \dfrac{1}{(n-1)h^p} \sum_{j \neq i} K_h(X_j - X_i)$.

We now turn to variance estimation. The leave-one-out construction of the $\widehat{T}_h$ gives that the asymptotic conditional variances $v_h^2$ and $v_{h,h_0}^2$ of $\widehat{T}_h$ and $\widehat{T}_h - \widehat{T}_{h_0}$ under $H_0$ are

$$
\begin{aligned}
(3.4) \qquad v_h^2 &= 2 \sum_{1 \leq i,j \leq n} w_{ij}^2(h) \sigma^2(X_i) \sigma^2(X_j), \\
v_{h,h_0}^2 &= 2 \sum_{1 \leq i,j \leq n} (w_{ij}(h) - w_{ij}(h_0))^2 \sigma^2(X_i) \sigma^2(X_j).
\end{aligned}
$$

For our main examples,

$$
v_{h_0}^2 \asymp_{\mathbb{P}} h_0^{-p} \quad \text{and} \quad v_{h,h_0}^2 \asymp_{\mathbb{P}} h^{-p} - h_0^{-p};
$$

see Proposition 2 in Section 6. Let $\sigma^2(\cdot)$ be a nonparametric estimator of $\widehat{\sigma}^2(\cdot)$ such that

$$
(3.5) \qquad \max_{1 \leq i \leq n} \left| \frac{\widehat{\sigma}_n^2(X_i)}{\sigma^2(X_i)} - 1 \right| = o_{\mathbb{P}}(1)
$$

for any equicontinuous sequence of regression functions. For instance, let

$$
\begin{aligned}
(3.6) \qquad \widehat{\sigma}_n^2(X_i) &= \frac{\sum_{j=1}^n Y_j^2 \mathbb{I}(|X_j - X_i| \leq b_n)}{\sum_{j=1}^n \mathbb{I}(|X_j - X_i| \leq b_n)} \\
&\quad - \left( \frac{\sum_{j=1}^n Y_j \mathbb{I}(|X_j - X_i| \leq b_n)}{\sum_{j=1}^n \mathbb{I}(|X_j - X_i| \leq b_n)} \right)^2,
\end{aligned}
$$

where $b_n$ is a bandwidth parameter chosen independently of $\mathcal{H}_n$ such that $n^{1-4/d'} b_n^p$ diverges; see Proposition 3 in Section 6. Consistent estimators of the variances in (3.4) are

$$
\begin{aligned}
\widehat{v}_{h_0}^2 &= 2 \sum_{1 \leq i,j \leq n} w_{ij}^2(h_0) \widehat{\sigma}_n^2(X_i) \widehat{\sigma}_n^2(X_j), \\
\widehat{v}_{h,h_0}^2 &= 2 \sum_{1 \leq i,j \leq n} (w_{ij}(h) - w_{ij}(h_0))^2 \widehat{\sigma}_n^2(X_i) \widehat{\sigma}_n^2(X_j).
\end{aligned}
$$

Finally, for the sake of parsimony, and following Horowitz and Spokoiny (2001), Lepski, Mammen and Spokoiny (1997) and Spokoiny (2001), the set $\mathcal{H}_n$ of admissible smoothing parameters is a geometric grid of $J_n + 1$ smoothing parameters,

$$
(3.7) \quad \mathcal{H}_n = \{h_j = h_0 a^{-j}, j = 0, \ldots, J_n\} \qquad \text{for some } a > 1, J_n \to +\infty.
$$

Note that $h_0$ can depend on an empirical measure of the dispersion of the $X_i$, as in Zhang (2003), and can converge to zero very slowly, say, as $1/\ln n$. We assume the following:

ASSUMPTION D. The i.i.d. $X_i \in [0,1]^p$ have a strictly positive continuous density over $[0,1]^p$.

ASSUMPTION M.    The function $\mu(x; \theta)$ is continuous with respect to $x$ in $[0, 1]^p$ and $\theta$ in $\Theta$, where $\Theta$ is a compact subset of $\mathbb{R}^d$. There is a constant $\dot{\mu}$ such that for all $\theta, \theta'$ in $\Theta$ and for all $x$ in $[0, 1]^p$, $|\mu(x; \theta) - \mu(x; \theta')| \leq \dot{\mu} |\theta - \theta'|$.

ASSUMPTION E.    The $\varepsilon_i$ are independent given $X_1, \ldots, X_n$. For each $i$, the distribution of $\varepsilon_i$, given the design, depends only on $X_i$, $\mathbb{E}[\varepsilon_i | X_i] = 0$ and $\mathrm{Var}[\varepsilon_i | X_i] = \sigma^2(X_i)$, where the unknown variance function $\sigma^2(\cdot)$ is continuous and bounded away from 0. For some $d' > \max(d, 4)$, $\mathbb{E}^{1/d'}[|\varepsilon_i|^{d'} | X_i] < C_1$ for all $i$.

ASSUMPTION W.    (i) For any $h$, the matrix $W_h$ is one from Example 1a, 1b or 2. (ii) The set $\mathcal{H}_n$ is as in (3.7) with $h_{J_n} \asymp (\ln n)^{C_2/p} n^{-2/(4\underline{s}+p)}$, for some $C_2 > 1$, with $\underline{s} = 5p/4$ in Example 1a and $\underline{s} = p/4$ in Examples 1b and 2. The number $a$ is an integer for Example 1b.

Under Assumption M, the value of the parameter $\theta$ may not be identified, as in mixture or multiple index models. The restriction on $h_{J_n}$, together with the definition of $\mathcal{H}_n$, implies that the number $J_n + 1$ of smoothing parameters is of order $\ln n$ at most. Assumption W(i), which considers specific nonparametric methods, will be relaxed in Section 5.1, allowing us, in particular, to consider a baseline statistic $\widehat{T}_{h_0}$ designed for specific parametric alternatives.

3.2. *Limit behavior of the test under the null hypothesis.*  The next theorem allows for a penalty sequence $\gamma_n$ of exact order $\sqrt{2 \ln \ln n}$, as $J_n$ is of order $\ln n$.

THEOREM 1.    *Consider a sequence* $\{\mu(\cdot, \theta_n), \theta_n \in \Theta\}_{n \geq 1}$ *in* $H_0$. *Let Assumptions* D, M, E *and* W *hold and assume that the variance estimator satisfies* (3.5). *If* $h_0 \to 0$ *and* $\gamma_n \to \infty$ *with*

$$(3.8) \qquad \gamma_n \geq (1 + \eta)\sqrt{2 \ln J_n} \qquad \textit{for some } \eta > 0,$$

*the test* (2.2) *has level* $\alpha$ *asymptotically given the design, that is,*

$$\mathbb{P}(\widehat{T}_{\widetilde{h}} \geq z_\alpha \widehat{v}_{h_0} | X_1, \ldots, X_n) \xrightarrow{\mathbb{P}} \alpha.$$

Theorem 1 is proved in two main steps. The first step consists in showing that

$$(3.9) \qquad \mathbb{P}(\widetilde{h} \neq h_0) = \mathbb{P}\left( \max_{h \in \mathcal{H}_n \setminus \{h_0\}} \frac{\widehat{T}_h - \widehat{T}_{h_0}}{\widehat{v}_{h,h_0}} > \gamma_n \right)$$

goes to zero. This is done by first proving that $(\widehat{T}_h - \widehat{T}_{h_0})/\widehat{v}_{h,h_0}$ asymptotically behaves at first-order as $\varepsilon'(W_h - W_{h_0})\varepsilon / v_{h,h_0}$ uniformly for $h$ in $\mathcal{H}_n \setminus \{h_0\}$, where $\varepsilon = [\varepsilon_1, \ldots, \varepsilon_n]'$, and second by bounding the distribution tails of $\max_{h \in \mathcal{H}_n \setminus \{h_0\}} \varepsilon'(W_h - W_{h_0})\varepsilon / v_{h,h_0}$. Then we show that the limit distribution of $\widehat{T}_{h_0} / \widehat{v}_{h_0}$ is that of $\varepsilon'W_{h_0}\varepsilon / v_{h_0}$, which converges to a standard normal when $h_0$ goes to 0.

As done by Horowitz and Spokoiny ([2001](#)), Theorem [1](#) imposes that $h_0$ asymptotically vanishes. This condition yields a pivotal limit distribution for our test statistic. As shown by Hart [([1997](#)), page 220] under stronger regularity conditions on the parametric model, considering a fixed $h_0$ generally yields a nonpivotal limit distribution because the estimation error $\mu(\cdot; \widehat{\theta}_n) - \mu(\cdot; \theta)$ cannot be neglected. Hart ([1997](#)) then recommends the use of a double bootstrap procedure to estimate the critical values of the test.

3.3. *Consistency of the test.* Theorem [2](#) below considers general alternatives with unknown smoothness. Theorem [3](#) considers Pitman local alternatives. For any real $s$, let $\lfloor s \rfloor$ be the lower integer part of $s$, that is, $\lfloor s \rfloor < s \le \lfloor s \rfloor + 1$. Let the Hölder class $C_p(L, s)$ be the set of maps $m(\cdot)$ from $[0,1]^p$ to $\mathbb{R}$ with

$$C_p(L, s) = \{m(\cdot); |m(x) - m(y)| \le L|x - y|^s \text{ for all } x, y \text{ in } [0,1]^p\}$$
$$\text{for } s \in (0, 1],$$
$$C_p(L, s) = \{m(\cdot); \text{ the } \lfloor s \rfloor\text{th partial derivatives of } m(\cdot) \text{ are in } C_p(L, s - \lfloor s \rfloor)\}$$
$$\text{for } s > 1.$$

THEOREM 2. *Consider a sequence of equicontinuous regression functions $\{m_n(\cdot)\}_{n \ge 1}$ such that for some unknown $s > \underline{s}$ and $L > 0$, $m_n(\cdot) - \mu(\cdot; \theta) \in C_p(L, s)$ for all $\theta$ in $\Theta$ and all $n$. Let Assumptions [D](#), [M](#), [E](#) and [W](#) hold. Assume that the variance estimator satisfies ([3.5](#)), that $1/(C_0 \ln n) \le h_0 \le C_0$ for some $C_0 > 0$ and that $\gamma_n \le n^\gamma$ for some $\gamma$ in $(0, 1)$. If*

$$(3.10) \quad \begin{aligned} &\min_{\theta \in \Theta} \left[ \frac{1}{n} \sum_{i=1}^n (m_n(X_i) - \mu(X_i; \theta))^2 \right]^{1/2} \\ &\ge (1 + o_{\mathbb{P}}(1)) \kappa_1 L^{p/(4s+p)} \left( \frac{\gamma_n \sup_{x \in [0,1]^p} \sigma^2(x)}{n} \right)^{2s/(4s+p)}, \end{aligned}$$

*the test ([2.2](#)) is consistent given the design, that is,*

$$\mathbb{P}(\widehat{T}_{\widehat{h}} \ge \widehat{v}_{h_0} z_\alpha | X_1, \ldots, X_n) \xrightarrow{\mathbb{P}} 1,$$

*provided $\kappa_1 = \kappa_1(s) > 0$ is large enough.*

The proof is based upon the power bound ([2.3](#)). From this inequality, the test is consistent if $\widehat{T}_h - z_\alpha \widehat{v}_{h_0} - \gamma_n \widehat{v}_{h,h_0}$ diverges in probability for a suitable choice of the smoothing parameter $h$ adapted to the unknown smoothness of the departure from the parametric model. Thus, combining several statistics in the procedure is crucial to detecting alternatives of unknown smoothness. A sketch of the proof is as follows. For a departure from the parametric model in $C_p(L, s)$, $\widehat{T}_h$ estimates $\min_{\theta \in \Theta} \sum_{i=1}^{n} (m_n(X_i) - \mu(X_i; \theta))^2$ up to a multiplicative constant with a bias of order $nL^2 h^{2s}$. The standard deviation of $\widehat{T}_h$ is of order $h^{-p/2}$ and the order of $\widehat{v}_{h_0} z_\alpha + \gamma_n \widehat{v}_{h,h_0}$ is $\gamma_n h^{-p/2} \sup_{x \in [0,1]^p} \sigma^2(x)$. Collecting the leading terms shows that $\widehat{T}_h - \widehat{v}_{h_0} z_\alpha - \gamma_n \widehat{v}_{h,h_0}$ diverges if

$$\min_{\theta \in \Theta} \left[ \frac{1}{n} \sum_{i=1}^{n} (m_n(X_i) - \mu(X_i; \theta))^2 \right]^{1/2}$$

is of larger order than

$$\left[ \frac{1}{n} \left( nL^2 h^{2s} + \gamma_n h^{-p/2} \sup_{x \in [0,1]^p} \sigma^2(x) \right) \right]^{1/2}.$$

Finding the minimum of this quantity with respect to $h$ gives the rate of (3.10). The rate of the optimal $h$ is $(\gamma_n \inf_{x \in [0,1]^p} \sigma^2(x) / L^2 n)^{2/(4s+p)}$. The parsimonious set $\mathcal{H}_n$ is rich enough to contain an $h$ of this order. Our proof can be easily modified to study the selection procedures considered in Hart (1997), which use $\gamma_n \widehat{v}_h^2$ in (2.1) instead of $\gamma_n \widehat{v}_{h,h_0}$. This would give the worst detection rate $(\gamma_n / n)^{s/(2s+p)}$.

For $\gamma_n$ of order $\sqrt{\ln \ln n}$, the smallest order compatible with Theorem 1, the test detects alternatives (3.10) with rate $(\sqrt{\ln \ln n}/n)^{2s/(4s+p)}$ for any $s > \underline{s}$. This rate is the optimal adaptive minimax one for the idealistic white noise model; see Spokoiny (1996). Horowitz and Spokoiny (2001) obtain the same rate for their kernel-based test but with minimal smoothness index $\underline{s} = \max(2, p/4)$, while we achieve $\underline{s} = p/4$ for our piecewise polynomial or kernel-based tests. The value $p/4$ is critical for the smoothness index $s$, as previously noted by Guerre and Lavergne (2002) and Baraud, Huet and Laurent (2003).

THEOREM 3. *Let $\theta_0$ be an inner point of $\Theta$ and consider a sequence of local alternatives $m_n(\cdot) = \mu(\cdot; \theta_0) + r_n \delta_n(\cdot)$, where $\{\delta_n(\cdot)\}_{n \geq 1}$ is an equicontinuous sequence from $C_p(L, s)$ for some unknown $s > \underline{s}$ and $L > 0$, with*

$$(3.11) \quad \frac{1}{n} \sum_{i=1}^{n} \delta_n^2(X_i) = 1 + o_{\mathbb{P}}(1) \quad and \quad \frac{1}{n} \sum_{i=1}^{n} \delta_n(X_i) \frac{\partial \mu(X_i; \theta_0)}{\partial \theta} = o_{\mathbb{P}}(1).$$

*Assume that for each $x$ in $[0,1]^p$, $\mu(x; \theta)$ is twice differentiable with respect to $\theta$ in $\Theta$ with second-order derivatives continuous in $x$ and $\theta$ and that, for*

*some $C_3 > 0$,*

$$(3.12) \quad \begin{aligned} &(C_3 + o_{\mathbb{P}}(1))|\theta - \theta'|^2 \\ &\leq \frac{1}{n}\sum_{i=1}^{n}(\mu(X_i;\theta) - \mu(X_i;\theta'))^2 \qquad \textit{for any } \theta, \theta' \textit{ in } \Theta. \end{aligned}$$

*Let Assumptions* D, M, E *and* W *hold and assume that the variance estimator satisfies* (3.5). *If* $h_0 \to 0$, $r_n \to 0$ *and* $\sqrt{nh_0^{p/2}}r_n \to \infty$, *the test is consistent given the design.*

The rate $r_n$ of Theorem 3 can be made arbitrarily close to $1/\sqrt{n}$ by a proper choice of $h_0$. This improves upon Horowitz and Spokoiny (2001), who obtain the rate $\sqrt{\ln\ln n}/\sqrt{n}$.

As stated in Lemma 5 of Section 6, conditions (3.11) and the identification condition (3.12) ensure that

$$(3.13) \quad \min_{\theta \in \Theta}\left[\frac{1}{n}\sum_{i=1}^{n}(m_n(X_i) - \mu(X_i;\theta))^2\right]^{1/2} = r_n - o_{\mathbb{P}}(r_n).$$

As the minimum of (3.13) is achieved for $\theta = \theta_0$ at first-order, $r_n\delta_n(\cdot)$ is asymptotically the departure from $\mu(\cdot;\theta_0)$. When $r_n$ converges to zero, this departure becomes smoother as it belongs to the smoothness class $C_p(Lr_n, s)$. This sharply contrasts with the departures from the parametric model in Theorem 2, which can be much more irregular. The proof of Theorem 3 follows from (2.4). The test is consistent as soon as $\widehat{T}_{h_0} - \widehat{v}_{h_0}z_\alpha$ diverges in probability. We show that $\widehat{T}_{h_0}$ is, up to a multiplicative constant, an estimate of $r_n^2\sum_{i=1}^{n}\delta_n^2(X_i)$ with a negligible bias and a standard deviation of order $h_0^{-p/2}$. As $\widehat{v}_{h_0}$ is of order $h_0^{-p/2}$, $\widehat{T}_{h_0} - \widehat{v}_{h_0}z_\alpha$ diverges to infinity as soon as $nr_n^2$ diverges faster than $h_0^{-p/2}$ as required.

## 4. Bootstrap implementation and small sample behavior.

4.1. *Bootstrap critical values.* The wild bootstrap, initially proposed by Wu (1986), is often used in smooth lack-of-fit tests to compute small sample critical values; see, for example, Härdle and Mammen (1993). Here we use a generalization of this method, the smooth conditional moments bootstrap introduced by Gozalo (1997). It consists of drawing $n$ i.i.d. random variables $\omega_i$ independently from the original sample with $\mathbb{E}\omega_i = 0$, $\mathbb{E}\omega_i^2 = 1$ and $\mathbb{E}|\omega_i|^{d'} < \infty$, and generating bootstrap observations of $Y_i$ as $Y_i^* = \mu(X_i, \widehat{\theta}_n) + \widehat{\sigma}_n(X_i)\omega_i, i = 1,\ldots,n$. A bootstrap test statistic $\widehat{T}_{\hat{h}^*}^*/\widehat{v}_{h_0}^*$ is built from the bootstrap sample, as was the original test statistic. When this scheme is repeated many times, the bootstrap critical value $z_{\alpha,n}^*$ at level $\alpha$ is the empirical $1 - \alpha$ quantile of the bootstrapped test statistics. This critical value is then compared to the initial test statistic. The following theorem establishes the first-order consistency of this procedure.

THEOREM 4. *Let* $Y_i = m_n(X_i) + \varepsilon_i$, $i = 1, \ldots, n$, *be the initial model, where* $\{m_n(\cdot)\}_{n \geq 1}$ *is any equicontinuous sequence of functions. Under the assumptions of Theorem* 1 *and for the variance estimator* $\hat{\sigma}_n^2(X_i)$ *of* (3.6),

$$\sup_{z \in \mathbb{R}} |\mathbb{P}(\widehat{T}_{\widetilde{h}^*}^* / \widehat{v}_{h_0}^* \leq z | X_1, Y_1, \ldots, X_n, Y_n) - \mathbb{P}(N(0,1) \leq z)| \xrightarrow{\mathbb{P}} 0.$$

4.2. *Small sample behavior.* We investigated the small sample behavior of our bootstrap test. We generated samples of 150 observations through the model

$$(4.1) \qquad Y = \theta_1 + \theta_2 X + r \cos(2\pi t X) + \varepsilon, \qquad r \in \left\{0, \sqrt{\tfrac{2}{3}}\right\}, t \in \{2, 5, 10\},$$

where $X$ is distributed as $U[-1, 1]$. The null hypothesis corresponds to $r = 0$, while under the alternatives $r^2 = 2/3$ and $\mathbb{E}[r^2 \cos^2(2\pi t X)]/\mathbb{E}\varepsilon^2 = 1/3$ for any integer $t$, a quite small signal-to-noise ratio. When $t$ increases, the deviation from the linear model becomes more oscillating and irregular, and then more difficult to detect.

To compute our test statistic, we used the regressogram method of Example 1b with half-binwidths in

$$\mathcal{H}_n = \{h_0 = 2^{-2}, h_1 = 2^{-3}, \ldots, h_5 = 2^{-7}\}.$$

The smallest binwidth thus defines 128 cells, which is sufficient for 150 observations. The $\gamma_n$ was set to $c\sqrt{2 \ln J_n}$, where $c = 1, 1.5, 2$. For each experiment we ran 5000 replications under the null and 1000 under the alternative. For each replication the bootstrap critical values were computed from 199 bootstrap samples. For $\omega_i$ we used the two-point distribution

$$\mathbb{P}\left(\omega_i = \frac{1 - \sqrt{5}}{2}\right) = \frac{5 + \sqrt{5}}{10}, \qquad \mathbb{P}\left(\omega_i = \frac{1 + \sqrt{5}}{2}\right) = \frac{5 - \sqrt{5}}{10},$$

which verifies the required conditions.

In a first stage we set $(\theta_1, \theta_2) = (0, 0)$ and performed a test for white noise, that is, $H_0 : m(\cdot) = 0$, with homoscedastic errors following a standard normal distribution (Table 1). We estimated the variance under homoscedasticity by

$$\hat{\sigma}_n^2 = \frac{1}{2(n-1)} \sum_{i=1}^{n-1} (Y_{(i+1)} - Y_{(i)})^2,$$

where $Y_{(i)}$ denote observations ordered according to the order of the $X_i$. This estimate is consistent under the null and the alternative; see Rice (1984).

In each cell of the tables, the first and second rows give empirical percentages of rejections at 2% and 5% nominal levels. We compare our test to (i) simple benchmark tests based on fixed bandwidths $h_0$ and $h_5$, to evaluate the effect of a data-driven bandwidth, (ii) the maximum test based on $\text{MAX} = \max_{h \in \mathcal{H}_n} \widehat{T}_h / \widehat{v}_h$, to evaluate the gain of our approach and (iii) a test based on $\widehat{T}_{\tilde{h}} / \widehat{v}_{\tilde{h}}$, to evaluate the effect of our standardization. For each test, we computed bootstrap critical values as for our test.

Under the null hypothesis, the bootstrap leads to accurate rejection probabilities for all tests. Under the considered alternatives, empirical power decreases for all tests when the frequency increases from $t = 2$ to $t = 10$. The data-driven tests always dominate the tests based on the fixed parameter $h_0$, which behaves poorly. For the low frequency alternatives, data-driven tests perform very well with power greater than 90% and 95% at a 2% and 5% nominal level, respectively, and there are no significant differences between them. For higher frequency alternatives, differences are significant. Our test has quite high power and rejects the null hypothesis at more than 85% and 60% at a 5% level when $t = 5$ and 10, respectively. It performs better than or as well as does the test based on $h_5$ designed for irregular alternatives, except for $c = 2$ and $t = 10$. It always dominates MAX with differences ranging from 7.1% to 18.3%, depending on the level. The test based on $\widehat{T}_{\tilde{h}} / \widehat{v}_{\tilde{h}}$ behaves as the MAX test. This suggests that the high performances of our test are mainly explained by our standardization choice, which is made possible by our selection procedure.

To check whether these conclusions are affected by the details of the experiments, we consider errors following a centered and standardized exponential

TABLE 1
*White noise model—Gaussian errors*

| | $\dfrac{\widehat{T}_{h_0}}{\widehat{v}_{h_0}}$ | $\dfrac{\widehat{T}_{h_{J_n}}}{\widehat{v}_{h_{J_n}}}$ | Max | $\dfrac{\widehat{T}_{\tilde{h}}}{\widehat{v}_{\tilde{h}}}$ | | | Our test | | |
|---|---|---|---|---|---|---|---|---|---|
| | | | | $c=1$ | $c=1.5$ | $c=2$ | $c=1$ | $c=1.5$ | $c=2$ |
| $H_0$ | 1.9 | 2.1 | 2.0 | 2.0 | 2.0 | 2.0 | 1.8 | 1.8 | 1.7 |
| | 5.3 | 5.1 | 4.2 | 4.3 | 4.2 | 4.4 | 4.4 | 4.3 | 4.4 |
| $t=2$ | 5.1 | 60.6 | 90.5 | 90.7 | 90.0 | 90.5 | 91.7 | 91.3 | 91.9 |
| | 9.0 | 72.5 | 96.0 | 96.3 | 95.9 | 96.2 | 95.4 | 95.7 | 97.3 |
| $t=5$ | 3.0 | 59.2 | 66.3 | 66.9 | 66.3 | 66.3 | 77.3 | 78.5 | 78.8 |
| | 7.7 | 73.3 | 79.2 | 79.8 | 79.4 | 79.5 | 88.7 | 88.5 | 87.8 |
| $t=10$ | 3.4 | 50.5 | 32.8 | 32.5 | 32.5 | 32.7 | 48.4 | 49.2 | 49.2 |
| | 7.0 | 66.0 | 49.3 | 50.2 | 49.3 | 48.8 | 65.6 | 65.5 | 59.9 |

Percentages of rejection at 2% and 5% nominal levels.

(Table 2), a standardized Student with five degrees of freedom (Table 3), a normal distribution with conditional variance $\sigma^2(X) = (1 + 3X^2)/3$ using our estimator (3.6) with $b_n = 1/8$ (Table 4) and a linear model with homoscedastic normal errors and $(\theta_1, \theta_2) = (1, 3)$ (Table 5). As results for $\widehat{T}_{\tilde{h}}/\widehat{v}_{\tilde{h}}$ are very similar to the ones for MAX, we do not report them. For exponential errors, there is a slight tendency to overrejection. It is likely that matching third-order moments in the bootstrap sample generation as proposed by Gozalo (1997) would lead to more accurate critical values. Heteroscedasticity does not adversely affect the behavior of the tests. For the linear model, there is some gain in power for the MAX test compared with Table 1, but differences with our test remain significant for the two high-frequency alternatives.

## 5. Extensions to general nonparametric methods and additive alternatives.

5.1. *General nonparametric methods.* We give here some general sufficient conditions ensuring the validity of our results. These conditions could be checked for other smoothing methods or other designs than the ones considered here. Indeed, different smoothing methods can be used for specification testing; see, for example, Chen (1994) for spline smoothing, Fan, Zhang and Zhang (2001) for local polynomials and Spokoiny (1996) for wavelets. Also, our conditions allow for various constructions of the quadratic forms $\widehat{T}_h$; see, for example, Dette (1999) and Härdle and Mammen (1993).

For an $n \times n$ matrix $W$, let $\mathrm{Sp}_n[W]$ be its spectral radius and $N_n^2[W] = \mathrm{Tr}[W'W] = \sum_{i,j} w_{ij}^2$. For $W$ symmetric, the former is its largest eigenvalue in absolute value and the latter is the sum of its squared eigenvalues.

TABLE 2
*White noise model—exponential errors*

| | $\dfrac{\widehat{T}_{h_0}}{\widehat{v}_{h_0}}$ | $\dfrac{\widehat{T}_{h,J_n}}{\widehat{v}_{h,J_n}}$ | **Max** | **Our test** | | |
|---|---|---|---|---|---|---|
| | | | | $c = 1$ | $c = 1.5$ | $c = 2$ |
| $H_0$ | 2.9 | 2.9 | 3.3 | 3.3 | 3.2 | 3.4 |
| | 6.1 | 6.2 | 6.7 | 6.3 | 5.9 | 6.5 |
| $t = 2$ | 4.5 | 65.4 | 91.9 | 92.2 | 92.4 | 92.6 |
| | 9.0 | 77.7 | 95.9 | 96.1 | 96.3 | 97.2 |
| $t = 5$ | 5.6 | 61.4 | 66.5 | 76.7 | 77.0 | 78.6 |
| | 9.6 | 71.7 | 78.9 | 86.1 | 87.0 | 86.0 |
| $t = 10$ | 3.6 | 50.6 | 35.4 | 51.3 | 52.8 | 53.7 |
| | 7.6 | 64.5 | 52.3 | 65.5 | 65.6 | 62.0 |

Percentages of rejection at 2% and 5% nominal levels.

TABLE 3
*White noise model—Student errors*

| | $\dfrac{\widehat{T}_{h_0}}{\widecheck{v}_{h_0}}$ | $\dfrac{\widehat{T}_{h_{J_n}}}{\widecheck{v}_{h_{J_n}}}$ | Max | **Our test** | | |
| | | | | $c=1$ | $c=1.5$ | $c=2$ |
|---|---|---|---|---|---|---|
| $H_0$ | 2.3 | 2.1 | 2.0 | 1.8 | 1.7 | 1.9 |
| | 5.0 | 4.8 | 4.4 | 4.5 | 4.3 | 4.4 |
| $t=2$ | 5.2 | 60.4 | 91.8 | 91.9 | 92.2 | 92.1 |
| | 9.2 | 73.3 | 95.7 | 95.5 | 95.8 | 96.2 |
| $t=5$ | 3.4 | 60.6 | 66.6 | 77.6 | 77.7 | 79.0 |
| | 8.4 | 74.6 | 79.3 | 88.2 | 88.2 | 86.9 |
| $t=10$ | 3.6 | 48.8 | 32.2 | 48.1 | 48.5 | 49.4 |
| | 7.8 | 65.1 | 48.1 | 63.1 | 64.2 | 60.0 |

Percentages of rejection at 2% and 5% nominal levels.

ASSUMPTION W0. Let $\mathcal{H}_n$ be as in (3.7) with $h_{J_n} \asymp (\ln n)^{C_2/p}/n^{2/(4\underline{s}+p)}$ for some $\underline{s} > 0$, $C_2 > 1$ and $h_0 \to 0$. The collection of $n \times n$ matrices $\{W_h, h \in \mathcal{H}_n\}$ is such that: (i) For all $h$, $W_h = [w_{ij}(h), 1 \le i, j \le n]$ depends only upon $X_1, \ldots, X_n$ and is real symmetric with $w_{ii}(h) = 0$ for all $i$. (ii) $\max_{h \in \mathcal{H}_n} \mathrm{Sp}_n[W_h] = O_{\mathbb{P}}(1)$. (iii) $N_n^2[W_h] \asymp_{\mathbb{P}} h^{-p}$ for all $h \in \mathcal{H}_n$ and uniformly in $h \in \mathcal{H}_n \setminus \{h_0\} N_n^2[W_h - W_{h_0}] \asymp_{\mathbb{P}} h^{-p} - h_0^{-p}$.

ASSUMPTION W1. Let $\mathcal{H}_n$, $\underline{s}$ and $h_{J_n}$ be as in Assumption W0. For any sequence $h_n = h_{j_n}$ from $\mathcal{H}_n$: (i) There are some symmetric positive semidefinite matrices $P_{h_n}$ with $\mathrm{Sp}_n[W_{h_n} - P_{h_n}] = o_{\mathbb{P}}(1)$. (ii) For any $s > \underline{s}$, there is

TABLE 4
*White noise model—heteroscedastic errors*

| | $\dfrac{\widehat{T}_{h_0}}{\widecheck{v}_{h_0}}$ | $\dfrac{\widehat{T}_{h_{J_n}}}{\widecheck{v}_{h_{J_n}}}$ | Max | **Our test** | | |
| | | | | $c=1$ | $c=1.5$ | $c=2$ |
|---|---|---|---|---|---|---|
| $H_0$ | 2.2 | 2.2 | 1.8 | 1.7 | 1.5 | 1.6 |
| | 5.1 | 5.0 | 4.7 | 4.2 | 4.1 | 4.2 |
| $t=2$ | 3.0 | 62.3 | 92.6 | 94.1 | 93.9 | 94.9 |
| | 5.9 | 76.3 | 98.0 | 97.9 | 98.4 | 98.7 |
| $t=5$ | 1.6 | 64.4 | 62.9 | 82.9 | 83.5 | 83.9 |
| | 4.2 | 78.9 | 81.9 | 91.9 | 92.8 | 91.6 |
| $t=10$ | 2.2 | 57.8 | 26.8 | 53.3 | 53.7 | 53.2 |
| | 5.6 | 72.8 | 50.3 | 69.5 | 71.3 | 63.5 |

Percentages of rejection at 2% and 5% nominal levels.

TABLE 5
*Linear model—Gaussian errors*

| | $\dfrac{\widehat{T}_{h_0}}{\widehat{v}_{h_0}}$ | $\dfrac{\widehat{T}_{h_{J_n}}}{\widehat{v}_{h_{J_n}}}$ | **Max** | **Our test** | | |
|---|---|---|---|---|---|---|
| | | | | $c = 1$ | $c = 1.5$ | $c = 2$ |
| $H_0$ | 2.3 | 2.1 | 1.9 | 1.9 | 2.0 | 2.0 |
| | 5.0 | 5.0 | 4.4 | 4.5 | 4.5 | 5.0 |
| $t = 2$ | 3.0 | 59.8 | 93.6 | 91.0 | 91.2 | 91.1 |
| | 6.3 | 71.7 | 96.7 | 95.5 | 95.6 | 96.8 |
| $t = 5$ | 2.7 | 58.2 | 73.2 | 77.7 | 77.9 | 78.5 |
| | 5.8 | 72.7 | 85.0 | 88.4 | 88.2 | 88.4 |
| $t = 10$ | 3.0 | 48.2 | 41.9 | 50.4 | 50.6 | 50.0 |
| | 7.0 | 64.4 | 58.8 | 66.0 | 66.2 | 61.8 |

Percentages of rejection at 2% and 5% nominal levels.

a set $\Pi_{s,n}$ of functions from $[0,1]^p$ to $\mathbb{R}$ such that for any $L > 0$ and any $\delta(\cdot)$ in $C_p(L,s)$, there is a $\pi(\cdot)$ in $\Pi_{s,n}$ with $\sup_{x \in [0,1]^p} |\delta(x) - \pi(x)| \le C_4 L h_n^s$ for some $C_4 = C_4(s) > 0$. (iii) Let $\Lambda_n^2 = \Lambda_n^2(s, h_n) = \inf_{\pi \in \Pi_{s,n}} \sum_{1 \le i,j \le n} \pi(X_i) p_{ij}(h_n) \pi(X_j) / \sum_{i=1}^n \pi(X_i)^2$, where $p_{ij}(h_n)$ is the generic element of $P_{h_n}$. For any $s > \underline{s}$, there is a constant $C_5 = C_5(s) > 0$ such that $\mathbb{P}(\Lambda_n > C_5) \to 1$.

Assumption W1 describes the approximation properties of the nonparametric method used to build the $W_h$ and allows us to extend a result of Ingster [(1993), page 253 and following]; see Lemma 6 in Section 6. The next proposition shows that our main examples satisfy Assumptions W0 and W1 under a regular i.i.d. random design.

PROPOSITION 1. *Assume that Assumption* D *holds, and let $\underline{s}$ be as in Assumption* W. *Then Examples* 1a, 1b *and* 2 *satisfy Assumptions* W0 *and* W1.

The next theorem extends our main results under Assumptions W0 and W1. In Section 6 we actually show Theorems 1–4 by proving Theorem 5 and Proposition 1.

THEOREM 5. *Theorems* 1 *and* 4 *hold under Assumption* W0 *in place of Assumptions* D *and* W. *Theorems* 2 *and* 3 *hold under Assumptions* W0 *and* W1 *in place of Assumptions* D *and* W.

5.2. *Additive alternatives.* Our general framework easily adapts to detection of specific alternatives. We focus here on additive nonparametric regressions $m(x) = m_1(x_1) + \cdots + m_p(x_p)$. The null hypothesis is

$$H_0 : m(\cdot) = \mu(\cdot; \theta) \qquad \text{for some } \theta \in \Theta,$$
$$\text{where } \mu(x; \theta) = \mu_1(x_1; \theta) + \cdots + \mu_p(x_p; \theta).$$

For ease of notation, we consider a modification of Example 1a where we remove cross-products of polynomial functions. Let $X_i = [X_{1i}, \ldots, X_{pi}]'$ and consider the $(p/h) \times n$ matrix $\Psi_h = [X_{1i}^k, \ldots, X_{pi}^k, i = 1, \ldots, n, k = 0, \ldots, 1/h]$. Let $W_h$ be the matrix obtained from $\Psi_h(\Psi_h' \Psi_h)^{-1} \Psi_h'$ by setting the diagonal entries to 0 and $\widehat{T}_h$ defined as in (3.3).

THEOREM 6. *Let the matrices $W_h$ be as above and $\mathcal{H}_n$ be as in* (3.7), *with $h_{J_n} \asymp (\ln n)^{C_6} / n^{1/3}$ for some $C_6 > 1$. Let Assumptions* D, E *and* M *hold. Consider a sequence of additive equicontinuous regression functions $\{m_n(\cdot)\}_{n \geq 1}$ and assume that the variance estimator satisfies* (3.5).

(i) *For $h_0$ and $\gamma_n$ as in Theorem* 1, *the test is asymptotically of level $\alpha$ given the design.*

(ii) *Assume that for some unknown $s > 5/4$ and $L > 0$, $m_n(\cdot) - \mu(\cdot; \theta)$ is in $C_p(L, s)$ for all $\theta$ in $\Theta$ and all $n$. For $h_0$ and $\gamma_n$ as in Theorem* 2 *and*

$$\min_{\theta \in \Theta} \left[ \frac{1}{n} \sum_{i=1}^{n} (m_n(X_i) - \mu(X_i; \theta))^2 \right]^{1/2}$$
$$\geq (1 + o_{\mathbb{P}}(1)) \kappa_2 L^{1/(4s+1)} \left( \frac{\gamma_n \sup_{x \in [0,1]} \sigma^2(x)}{n} \right)^{2s/(4s+1)},$$

*the test is consistent given the design provided $\kappa_2 = \kappa_2(s)$ is large enough.*

Proof of Theorem 6 repeats the proofs of Theorems 1 and 2 with $v_{h,h_0}^2$ of order $(h^{-1} - h_0^{-1})$ instead of $(h^{-p} - h_0^{-p})$ and is therefore omitted. One can also show consistency of the test against Pitman additive alternatives that approach the parametric model at rate $o(1/\sqrt{n h_0^{1/2}})$. The bootstrap procedure described in Section 4.1 also remains valid.

**6. Proofs.** This section is organized as follows. In Section 6.1 we study the quadratic forms $\varepsilon'(W_h - W_{h_0})\varepsilon$ and $\varepsilon' W_h \varepsilon$ under $H_0$. Section 6.2 recalls some results related to variance estimation. In Section 6.3 we gather preliminary results on the parametric estimation error $m_n(\cdot) - \mu(\cdot; \widehat{\theta}_n)$. In Sections 6.4 and 6.5 we establish Theorems 1 and 4 under Assumption W0. In Sections 6.6 and 6.7 we establish Theorems 2 and 3 under Assumptions W0 and W1. Thus, Theorem 5 is a direct consequence of Sections 6.4–6.7. Section 6.8 deals with Proposition 1.

We denote $Y = [Y_1, \ldots, Y_n]'$ and $\varepsilon = [\varepsilon_1, \ldots, \varepsilon_n]'$. For any $\delta(\cdot)$ from $\mathbb{R}^p$ to $\mathbb{R}$, $\delta = \delta(X) = [\delta(X_1), \ldots, \delta(X_n)]'$ and $D_n(\delta)$ is the $n \times n$ diagonal matrix with entries $\delta(X_i)$. Let $\| \cdot \|_n^2$ and $(\cdot, \cdot)_n$ be the Euclidean norm and inner product on $\mathbb{R}^n$ divided by $n$, respectively, that is,

$$\|\delta\|_n^2 = \|\delta(X)\|_n^2 = \frac{1}{n} \sum_{i=1}^n \delta^2(X_i)$$

and

$$(\varepsilon, \delta)_n = (\varepsilon, \delta(X))_n = \frac{1}{n} \sum_{i=1}^n \varepsilon_i \delta(X_i).$$

This gives $\mathrm{Sp}_n[W] = \max_{\|u\|_n = 1} \|Wu\|_n = \max_{\|u\|_n = 1} |u'Wu|/n$ for a symmetric $W$. Recall that $\mathrm{Sp}_n[AB] \leq \mathrm{Sp}_n[A]\mathrm{Sp}_n[B]$. Let $\theta_n = \theta_{n,m}$ be such that

$$(6.1) \qquad \min_{\theta \in \Theta} \|m(X) - \mu(X; \theta)\|_n = \|m(X) - \mu(X; \theta_n)\|_n.$$

We use the notation $\mathbb{P}_n(A)$ for $\mathbb{P}(A|X_1, \ldots, X_n)$, $\mathbb{E}_n[\cdot]$ and $\mathrm{Var}_n[\cdot]$ being the associated conditional mean and variance. In what follows, $C$ and $C'$ are positive constants that may vary from line to line. An absolute constant depends neither on the design nor on the distribution of the $\varepsilon_i$ given the design.

6.1. *Study of quadratic forms.* The proof of Lemma 1 is omitted.

LEMMA 1. *Let $W$ be an $n \times n$ symmetric matrix with zeros on the diagonal. Under Assumption* E, $\mathbb{E}_n[\varepsilon'W\varepsilon] = 0$ *and* $\mathrm{Var}_n[\varepsilon'W\varepsilon] = 2\sum_{1 \leq i,j \leq n} w_{ij}^2 \sigma^2(X_i) \times \sigma^2(X_j) = 2N_n^2[D_n(\sigma)WD_n(\sigma)] \asymp N_n^2[W]$.

LEMMA 2. *Let $\underline{\sigma} = \inf_{x \in [0,1]^p} \sigma(x) > 0$, $\overline{\sigma} = \sup_{x \in [0,1]^p} \sigma(x) < \infty$ and $\nu \in (0, 1/2)$. Under Assumption* E, *there is an absolute constant $C = C_\nu > 0$ such that:*

(i) *If* $(\overline{\sigma}^4 \mathrm{Sp}_n^2[W_h])/(\underline{\sigma}^4 N_n^2[W_h]) \leq \nu$,

$$\sup_{z \in \mathbb{R}} |\mathbb{P}_n(\varepsilon'W_h\varepsilon \leq v_h z) - \mathbb{P}(N(0,1) \leq z)| \leq C\left(\frac{\overline{\sigma}\,\mathrm{Sp}_n[W_h]}{\underline{\sigma}\,N_n[W_h]}\right)^{1/4}.$$

(ii) *For all $h \in \mathcal{H}_n \setminus \{h_0\}$ and any $z > 0$, if $(\overline{\sigma}^4 \mathrm{Sp}_n^2[W_h - W_{h_0}])/(\underline{\sigma}^4 N_n^2[W_h - W_{h_0}]) < \nu$,*

$$\mathbb{P}_n\left(\left|\frac{\varepsilon'(W_h - W_{h_0})\varepsilon}{v_{h,h_0}}\right| \geq z\right) \leq \frac{\sqrt{2}}{\sqrt{\pi}z}\exp\left(-\frac{z^2}{2}\right) + C\left(\frac{\overline{\sigma}\,\mathrm{Sp}_n[W_h - W_{h_0}]}{\underline{\sigma}\,N_n[W_h - W_{h_0}]}\right)^{1/4}.$$

PROOF. Let $\widetilde{\varepsilon} = D_n^{-1}(\sigma)\varepsilon$, so that $\mathbb{E}_n[\widetilde{\varepsilon}_i] = 0$ and $\mathrm{Var}_n[\widetilde{\varepsilon}_i] = 1$ for all $i$, and let $W = [w_{ij}]_{1 \leq i, j \leq n}$ be $D_n(\sigma)W_h D_n(\sigma)$ or $D_n(\sigma)(W_h - W_{h_0})D_n(\sigma)$, so that for $v^2 = N_n^2[W] = \sum_{1 \leq i, j \leq n} w_{ij}^2$, $\widetilde{\varepsilon}'W\widetilde{\varepsilon}/v$ is $\varepsilon' W_h \varepsilon/v_h$ or $\varepsilon'(W_h - W_{h_0})\varepsilon/v_{h,h_0}$, respectively. Let $\widetilde{\lambda}_1, \ldots, \widetilde{\lambda}_n$ be the real eigenvalues of $W$,

$$\mathcal{L}_n = \frac{1}{v^3}\left[6\sum_{i=1}^n\left(\sum_{j=1}^n w_{ij}^2\right)^{3/2} + 36\sum_{i=1}^n\sum_{j=1}^n |w_{ij}|^3\right] \quad \text{and} \quad \Delta_n = \frac{1}{v^4}\sum_{i=1}^n \lambda_i^4.$$

Consider a vector $g$ of $n$ independent $N(0,1)$ variables, independent of the $X_i$. Theorem 3 of Rotar' and Shervashidze ([1985](#)) says that there is an absolute constant $C > 0$ such that

$$\sup_{z \in \mathbb{R}}\left|\mathbb{P}_n\left(\frac{\widetilde{\varepsilon}'W\widetilde{\varepsilon}}{v} \leq z\right) - \mathbb{P}_n\left(\frac{g'Wg}{v} \leq z\right)\right|$$
$$\leq C[1 - \ln(1 - 2\Delta_n)]^{3/4}\mathcal{L}_n^{1/4} \quad \text{if } \Delta_n < 1/2.$$

Let $\{b_i \in \mathbb{R}^n\}_{1 \leq i \leq n}$ be an orthonormal system of eigenvectors of $W$ associated with the eigenvalues $\lambda_i$. As $\mathbb{E}_n[g'Wg] = 0$ by Lemma [1](#), $g'Wg = \sum_{i=1}^n \lambda_i(b_i'g)^2 = \sum_{i=1}^n \lambda_i[(b_i'g)^2 - \mathbb{E}_n[(b_i'g)^2]]$. Hence, $g'Wg$ has the same conditional distribution as $\sum_{i=1}^n \lambda_i\zeta_i$, where the $\zeta_i$ are centered Chi-squared variables with one degree of freedom, independent among themselves and of the $X_i$. The Berry–Esseen bound of Chow and Teicher [([1988](#)), Theorem 3, page 304] yields that there is an absolute constant $C > 0$ such that

$$\sup_{z \in \mathbb{R}}\left|\mathbb{P}_n\left(\frac{g'Wg}{v} \leq z\right) - \mathbb{P}(N(0,1) \leq z)\right| \leq C\frac{\sum_{i=1}^n |\lambda_i|^3}{v^3}.$$

The two above inequalities together imply that if $\Delta_n < 1/2$,

$$(6.2) \quad \begin{aligned}\sup_{z \in \mathbb{R}}&\left|\mathbb{P}_n\left(\frac{\widetilde{\varepsilon}'W\widetilde{\varepsilon}}{v} \leq z\right) - \mathbb{P}(N(0,1) \leq z)\right| \\ &\leq C\left[(1 - \ln(1 - 2\Delta_n))^{3/4}\mathcal{L}_n^{1/4} + \frac{\sum_{i=1}^n |\lambda_i|^3}{v^3}\right].\end{aligned}$$

Let $\{e_i, i = 1, \ldots, n\}$ be the canonical basis of $\mathbb{R}^n$, so that $\|e_i\|_n = 1/\sqrt{n}$. Then

$$\sum_{i=1}^n\left(\sum_{j=1}^n w_{ij}^2\right)^{3/2} = \sum_{i=1}^n \frac{\|We_i\|_n}{\|e_i\|_n}n\|We_i\|_n^2$$
$$\leq \mathrm{Sp}_n[W] \times \sum_{1 \leq i, j \leq n} w_{ij}^2 = \mathrm{Sp}_n[W]N_n^2[W],$$

$$\sum_{1 \leq i, j \leq n} |w_{ij}|^3 = \sum_{1 \leq i, j \leq n} w_{ij}^2\frac{|(e_i, We_j)_n|}{\|e_i\|_n\|e_j\|_n}$$
$$\leq \sum_{1 \leq i, j \leq n} w_{ij}^2\frac{\|We_j\|_n}{\|e_j\|_n} \leq \mathrm{Sp}_n[W]N_n^2[W].$$

Hence, using $v^2 = \sum_{i=1}^n \lambda_i^2 = N_n^2[W]$ and $|\lambda_i| \leq \mathrm{Sp}_n[W]$ for all $i$, we obtain

$$\Delta_n \leq \frac{\mathrm{Sp}_n^2[W]}{N_n^2[W]}, \qquad \mathcal{L}_n \leq 42 \frac{\mathrm{Sp}_n[W]}{N_n[W]}$$

and

$$\sum_{i=1}^n \frac{|\lambda_i|^3}{v^3} \leq \frac{\mathrm{Sp}_n[W]}{N_n[W]} \leq \left( \frac{\mathrm{Sp}_n[W]}{N_n[W]} \right)^{1/4},$$

since $\mathrm{Sp}_n[W]/N_n[W] \leq 1$ for any symmetric $W$. The above inequalities and (6.2) give

$$(6.3) \qquad \sup_{z \in \mathbb{R}} \left| \mathbb{P}_n \left( \frac{\widetilde{\varepsilon}' W \widetilde{\varepsilon}}{v} \leq z \right) - \mathbb{P}(N(0,1) \leq z) \right| \leq C \left( \frac{\mathrm{Sp}_n[W]}{N_n[W]} \right)^{1/4},$$

provided $(\mathrm{Sp}_n[W]/N_n[W])^2 \leq \nu$, for an absolute constant $C = C_\nu > 0$.

Part (i) follows by setting $W = D_n(\sigma) W_h D_n(\sigma)$ in (6.3) and noting that

$$\left( \frac{\mathrm{Sp}_n[W]}{N_n[W]} \right)^2 \leq \left( \frac{\overline{\sigma}}{\underline{\sigma}} \right)^4 \left( \frac{\mathrm{Sp}_n[W_h]}{N_n[W_h]} \right)^2 \leq \nu < 1/2.$$

Part (ii) follows from (6.3) with $W = D_n(\sigma)(W_h - W_{h_0}) D_n(\sigma)$ and Mills' ratio inequality. $\square$

6.2. *Variance estimation.* The following results are proven in Guerre and Lavergne (2003).

PROPOSITION 2. *Under Assumptions* D *and* W, *$v_{h_0}^2 \asymp_{\mathbb{P}} h_0^{-p}$ and uniformly in $h \in \mathcal{H}_n \setminus \{h_0\}$ $v_{h,h_0}^2 \asymp_{\mathbb{P}} h^{-p} - h_0^{-p}$.*

PROPOSITION 3. *Let $\{m_n(\cdot)\}_{n \geq 1}$ be an equicontinuous sequence of regression functions.*

(i) *Under Assumptions* D *and* E, *if $b_n \to 0$ and $n^{1-4/d'} b_n^p \to \infty$, then* (3.5) *holds.*

(ii) *Let $\{W_h, h \in \mathcal{H}_n\}$ be any collection of nonzero $n \times n$ symmetric matrices with zeros on the diagonal. Under* (3.5), *$\frac{\widehat{v}_{h_0}^2}{v_{h_0}^2} \xrightarrow{\mathbb{P}} 1$ and $\max_{h \in \mathcal{H}_n \setminus \{h_0\}} |\frac{\widehat{v}_{h,h_0}^2}{v_{h,h_0}^2} - 1| = o_{\mathbb{P}}(1)$.*

6.3. *The parametric estimation error.*

LEMMA 3. *Let $W$ be an $n \times n$ symmetric matrix depending upon $X_1, \ldots, X_n$, $\theta_n$ be as in* (6.1) *and $B_n(R) = \{\theta \in \Theta; \frac{1}{n} \sum_{i=1}^n (\mu(X_i; \theta) - \mu(X_i; \theta_n))^2 \leq R^2\}$.*

*Under Assumptions* E *and* M, *there is an absolute constant* $C = C_{d'} > 0$ *such that, for any* $m_n(\cdot)$, *any* $n$ *and any* $R > 0$,

$$\mathbb{E}_n \left[ \sup_{\theta \in B_n(R)} |\sqrt{n}(W(\mu(X;\theta) - \mu(X;\theta_n)), \varepsilon)_n| \right]$$

$$\leq C\dot{\mu} \mathrm{Sp}_n[W] R \max_{1 \leq i \leq n} \mathbb{E}_n^{1/d'}[|\varepsilon_i|^{d'}].$$

PROOF. Without loss of generality, we can assume that $\max_{1 \leq i \leq n} \mathbb{E}^{1/d'}[|\varepsilon_i|^{d'} | X_i] = \dot{\mu} = \mathrm{Sp}_n[W] = 1$. Let $\delta_W(\cdot; \theta) = W(\mu(\cdot; \theta) - \mu(\cdot; \theta_n))$. The Marcinkie- wicz–Zygmund inequality, see Chow and Teicher ([1988](#)), yields, under Assumption E and for any $\theta, \theta'$ in $\Theta$, that there is an absolute constant $C$ such that

$$\mathbb{E}_n^{1/d'} \left| \frac{1}{\sqrt{n}} \sum_{i=1}^{n} (\delta_W(X_i; \theta) - \delta_W(X_i; \theta'))\varepsilon_i \right|^{d'}$$

$$\leq C \left[ \frac{1}{n} \sum_{i=1}^{n} (\delta_W(X_i; \theta) - \delta_W(X_i; \theta'))^2 \mathbb{E}_n^{2/d'} |\varepsilon_i|^{d'} \right]^{1/2}$$

$$\leq C \|W(\mu(X;\theta) - \mu(X;\theta'))\|_n \leq C \|\mu(X;\theta) - \mu(X;\theta')\|_n.$$

Let $\mathcal{N}_n(t, R)$ be the smallest number of $\|\mu(X;\theta) - \mu(X;\theta')\|_n$-balls of radius $t$ covering $B_n(R)$. It follows from van der Vaart [([1998](#)), Example 19.7] and Assumption M that, for some absolute constant $C' > 0$, $\mathcal{N}_n(t, R) \leq C'(R/t)^d$. The Hölder inequality and Corollary 2.2.5 from van der Vaart and Wellner ([1996](#)) give, as $d/d' < 1$,

$$\mathbb{E}_n \sup_{\theta \in B_n(R)} \left| \frac{1}{\sqrt{n}} \sum_{i=1}^{n} \delta_W(X_i; \theta)\varepsilon_i \right| \leq \mathbb{E}_n^{1/d'} \sup_{\theta \in B_n(R)} \left| \frac{1}{\sqrt{n}} \sum_{i=1}^{n} \delta_W(X_i; \theta)\varepsilon_i \right|^{d'}$$

$$\leq C' \int_0^R \left( \frac{R}{t} \right)^{d/d'} dt = C_{d'} R. \qquad \square$$

LEMMA 4. *Under Assumptions* E *and* M, *there is an absolute constant* $C = C_{d'} > 0$, *such that, for any* $\rho$ *large enough, any* $m_n(\cdot)$ *and any* $n$,

$$\mathbb{P}_n \left[ \|m_n(X) - \mu(X; \hat{\theta}_n)\|_n > \sqrt{3} \|m_n(X) - \mu(X; \theta_n)\|_n + \frac{\sqrt{2}\rho}{\sqrt{n}} \right]$$

$$\leq \frac{C \max_{1 \leq i \leq n} \mathbb{E}_n^{1/d'}[|\varepsilon_i|^{d'}]}{\rho}.$$

PROOF. The definition (3.1) of $\widehat{\theta}_n$ yields, see van de Geer (2000),

$$(6.4)\quad\begin{aligned}&\|m_n(X) - \mu(X;\widehat{\theta}_n)\|_n^2\\&\quad\leq 2(\mu(X;\widehat{\theta}_n) - \mu(X;\theta_n), \varepsilon)_n + \|m_n(X) - \mu(X;\theta_n)\|_n^2,\\&\|\mu(X;\widehat{\theta}_n) - \mu(X;\theta_n)\|_n^2\\&\quad\leq 4(\mu(X;\widehat{\theta}_n) - \mu(X;\theta_n), \varepsilon)_n + 4\|m_n(X) - \mu(X;\theta_n)\|_n^2.\end{aligned}$$

Consider a fixed $r > 1$ and any $\rho \geq r$. Let $\mathcal{E}_n = \{\|m_n(X) - \mu(X;\theta_n)\|_n^2 < (\mu(X;\widehat{\theta}_n) - \mu(X;\theta_n), \varepsilon)_n\}$, so that on the complement of this event $\|m_n(X) - \mu(X;\widehat{\theta}_n)\|_n \leq \sqrt{3}\|m_n(X) - \mu(X;\theta_n)\|_n$ by (6.4). Lemma 4 follows by bounding

$$\mathbb{P}_n\left(\left(\sqrt{3}\|m_n(X) - \mu(X;\theta_n)\|_n + \frac{\sqrt{2}r^J}{\sqrt{n}}\right)^2 \leq \|m_n(X) - \mu(X;\widehat{\theta}_n)\|_n^2 \text{ and } \mathcal{E}_n\right)$$

$$\leq \mathbb{P}_n\left(2\|m_n(X) - \mu(X;\theta_n)\|_n^2 + \frac{2r^{2J}}{n}\right.$$

$$\left.\leq 2\|m_n(X) - \mu(X;\theta_n)\|_n^2 + 2\|\mu(X;\theta_n) - \mu(X;\widehat{\theta}_n)\|_n^2 \text{ and } \mathcal{E}_n\right)$$

$$= \mathbb{P}_n\left(\frac{r^{2J}}{n} \leq \|\mu(X;\widehat{\theta}_n) - \mu(X;\theta_n)\|_n^2 \text{ and } \mathcal{E}_n\right).$$

Let $S_j = S_{j,n} = \{\theta \in \Theta; r^j/\sqrt{n} \leq \|\mu(X;\theta) - \mu(X;\theta_n)\|_n < r^{j+1}/\sqrt{n}\} \subset B_n(r^{j+1}/\sqrt{n})$ with $B_n(\cdot)$ as in Lemma 3. Then (6.4), the definition of $\mathcal{E}_n$, the Markov inequality and Lemma 3 with $W = \mathrm{Id}_n$ yield

$$\mathbb{P}_n\left(\frac{r^{2J}}{n} \leq \|\mu(X;\widehat{\theta}_n) - \mu(X;\theta_n)\|_n^2 \text{ and } \mathcal{E}_n\right)$$

$$\leq \sum_{j=J}^{+\infty} \mathbb{P}_n\left(\widehat{\theta}_n \in S_j \text{ and } \frac{r^{2j}}{8n} \leq (\mu(X;\widehat{\theta}_n) - \mu(X;\theta_n), \varepsilon)_n\right)$$

$$\leq \sum_{j=J}^{+\infty} \mathbb{P}_n\left(\frac{r^{2j}}{8\sqrt{n}} \leq \sup_{\theta \in B_n(r^{j+1}/\sqrt{n})} |\sqrt{n}(\mu(X;\theta) - \mu(X;\theta_n), \varepsilon)_n|\right)$$

$$\leq \sum_{j=J}^{+\infty} \frac{8\sqrt{n}}{r^{2j}} \mathbb{E}_n\left[\sup_{\theta \in B_n(r^{j+1}/\sqrt{n})} |\sqrt{n}(\mu(X;\theta) - \mu(X;\theta_n), \varepsilon)_n|\right]$$

$$\leq C \max_{1 \leq i \leq n} \mathbb{E}_n^{1/d'}[|\varepsilon_i|^{d'}] \sum_{j=J}^{+\infty} \frac{r^{j+1}\sqrt{n}}{r^{2j}\sqrt{n}}$$

$$= \frac{r^2}{r-1} \frac{C \max_{1 \leq i \leq n} \mathbb{E}_n^{1/d'}[|\varepsilon_i|^{d'}]}{r^J}. \qquad \square$$

Lemma 5 is proven in Guerre and Lavergne (2003).

LEMMA 5. *Consider the local alternatives of Theorem 3 and let the conditions of Theorem 3 on $\mu(\cdot; \cdot)$ hold. Under Assumptions E and M and if $\lim_{n \to +\infty} \sqrt{n} r_n = +\infty$,*

$$\|m_n(X) - \mu(X; \theta_n)\|_n = r_n - o_{\mathbb{P}}(r_n) \quad and \quad \|\mu(X; \widehat{\theta}_n) - \mu(X; \theta_0)\|_n = o_{\mathbb{P}}(r_n).$$

PROPOSITION 4. *Under Assumptions E, M and W0(ii), if $h_0 \to 0$, then, for any $\{m_n(\cdot)\}_{n \geq 1} \subset H_0$,*

$$\max_{h \in \mathcal{H}_n \setminus \{h_0\}} \left| \frac{\widehat{T}_h - \widehat{T}_{h_0} - \varepsilon'(W_h - W_{h_0})\varepsilon}{(h^{-p} - h_0^{-p})^{1/2}} \right| = o_{\mathbb{P}}(1), \qquad h_0^{p/2}(\widehat{T}_{h_0} - \varepsilon' W_{h_0}\varepsilon) = o_{\mathbb{P}}(1).$$

*Let $h_n \in \mathcal{H}_n$ be an arbitrary sequence of smoothing parameters. Then under $H_0$ or $H_1$,*

$$(m_n(X) - \mu(X, \widehat{\theta}_n))' W_h \varepsilon = O_{\mathbb{P}}(1)[\sqrt{n}\|m_n(X) - \mu(X, \theta_n)\|_n + 1].$$

PROOF. We have

$$(6.5) \qquad \begin{aligned} \widehat{T}_h &= (m_n(X) - \mu(X; \widehat{\theta}_n))' W_h(m_n(X) - \mu(X; \widehat{\theta}_n)) \\ &\quad + 2(m_n(X) - \mu(X; \widehat{\theta}_n))' W_h \varepsilon + \varepsilon' W_h \varepsilon. \end{aligned}$$

The Cauchy–Schwarz inequality, Assumptions E and W0(ii) and Lemma 4 yield uniformly in $h \in \mathcal{H}_n$,

$$\begin{aligned} &|(m_n(X) - \mu(X; \widehat{\theta}_n))' W_h(m_n(X) - \mu(X; \widehat{\theta}_n))| \\ &\qquad \leq n \max_{h \in \mathcal{H}_n} \mathrm{Sp}_n[W_h]\|m_n(X) - \mu(X; \widehat{\theta}_n)\|_n^2 \\ &\qquad = O_{\mathbb{P}}[(1 + \sqrt{n}\|m_n(X) - \mu(X; \theta_n)\|_n)^2] = O_{\mathbb{P}}(1) \end{aligned}$$

under $H_0$, as $\|m_n(X) - \mu(X; \theta_n)\|_n = 0$. Since for any $h \in \mathcal{H}_n$, $h^{-p} - h_0^{-p} \geq h_1^{-p} - h_0^{-p} = h_0^{-p}(a^p - 1) \to +\infty$, we obtain that, under $H_0$,

$$\begin{aligned} &\max_{h \in \mathcal{H}_n \setminus \{h_0\}} \left| \frac{(m_n(X) - \mu(X; \widehat{\theta}_n))'(W_h - W_{h_0})(m_n(X) - \mu(X; \widehat{\theta}_n))}{(h^{-p} - h_0^{-p})^{1/2}} \right| = o_{\mathbb{P}}(1), \\ &\qquad\qquad h_0^{p/2}(m_n(X) - \mu(X; \widehat{\theta}_n))' W_{h_0}(m_n(X) - \mu(X; \widehat{\theta}_n)) = o_{\mathbb{P}}(1). \end{aligned}$$
$$(6.6)$$

Since $\|\mu(X; \widehat{\theta}_n) - \mu(X; \theta_n)\|_n \leq \|\mu(X; \widehat{\theta}_n) - m_n(X)\|_n + \|m_n(X) - \mu(X; \theta_n)\|_n$, Lemma 4 and Assumption E yield $\mathbb{P}_n(\widehat{\theta}_n \notin B_{\rho,n}) \leq C/\rho$ for any $\rho$ large enough, any $m_n(\cdot)$ and any $n$, where

$$\begin{aligned} B_{\rho,n} = \Big\{ \theta \in \Theta; \\ \|\mu(X; \theta) - \mu(X; \theta_n)\|_n \leq (\sqrt{3} + 1)\|m_n(X) - \mu(X; \theta_n)\|_n + \frac{\sqrt{2}\rho}{\sqrt{n}} \Big\}. \end{aligned}$$

Lemma 3 yields

$$
\begin{aligned}
(6.7)\qquad & \mathbb{E}_n\Big[\sup_{\theta\in B_{\rho,n}}|(\mu(X,\theta)-\mu(X;\theta_n))'W\varepsilon|\Big] \\
& \le C\rho\mathrm{Sp}_n[W](\sqrt{n}\|m_n(X)-\mu(X;\theta_n)\|_n+1).
\end{aligned}
$$

Taking $W=W_{h_0}$ and using the Markov inequality, (6.5), (6.6), $m_n(X)-\mu(X;\theta_n)=0$, Assumption W0(ii) and $h_0\to0$ then show that $h_0^{p/2}(\widehat{T}_{h_0}-\varepsilon'W_{h_0}\varepsilon)=o_{\mathbb{P}}(1)$ under $H_0$. Taking $W=W_h-W_{h_0}$ in (6.7) and using $h=h_0a^{-j}$ for some $j=0,\dots,J_n$ yields, under $H_0$,

$$
\begin{aligned}
& \mathbb{P}_n\Big(\max_{h\in\mathcal{H}_n\setminus\{h_0\}}\Big|\frac{(\mu(X,\widehat{\theta}_n)-\mu(X;\theta_n))'(W_h-W_{h_0})\varepsilon}{(h^{-p}-h_0^{-p})^{1/2}}\Big|\ge\epsilon\Big) \\
& \quad\le \mathbb{P}_n(\widehat{\theta}_n\notin B_{\rho,n}) \\
& \qquad + \frac{1}{\epsilon}\sum_{h\in\mathcal{H}_n\setminus\{h_0\}}\mathbb{E}_n\sup_{\theta\in B_{\rho,n}}\Big|\frac{(\mu(X,\theta)-\mu(X;\theta_n))'(W_h-W_{h_0})\varepsilon}{(h^{-p}-h_0^{-p})^{1/2}}\Big| \\
& \quad\le \frac{C}{\rho}+\frac{\rho}{\epsilon}O_{\mathbb{P}}(h_0^{p/2})\sum_{j=1}^{\infty}\frac{1}{(a^{pj}-1)^{1/2}}=\frac{C}{\rho}+\frac{\rho}{\epsilon}O_{\mathbb{P}}(h_0^{p/2}),
\end{aligned}
$$

for all $\epsilon>0$. The last result follows from (6.7) with $W=W_h$ and

$$
\mathbb{E}_n[((m_n(X)-\mu(X;\theta_n))'W_h\varepsilon)^2]\le n\mathrm{Sp}_n^2(W_h)\overline{\sigma}^2\|m_n(X)-\mu(X;\theta_n)\|_n^2.\quad\square
$$

6.4. *Proof of Theorem* 1 *under Assumption* W0. Under Assumptions W0(iii) and E, $v_{h,h_0}\asymp N_n[W_h-W_{h_0}]\asymp_{\mathbb{P}}(h^{-p}-h_0^{-p})^{1/2}$ uniformly in $h\in\mathcal{H}_n\setminus\{h_0\}$; see Lemma 1. Therefore, Propositions 3(ii) and 4 yield

$$
\max_{h\in\mathcal{H}_n\setminus\{h_0\}}\Big|\frac{\widehat{T}_h-\widehat{T}_{h_0}}{\widehat{v}_{h,h_0}}\Big|=(1+o_{\mathbb{P}}(1))\times\max_{h\in\mathcal{H}_n\setminus\{h_0\}}\Big|\frac{\varepsilon'(W_h-W_{h_0})\varepsilon}{v_{h,h_0}}\Big|+o_{\mathbb{P}}(1).
$$

Let $\eta$ be as in (3.8). Observe that

$$
\begin{aligned}
\mathbb{P}_n(\widetilde{h}\neq h_0)&\le\mathbb{P}_n\Big(\max_{h\in\mathcal{H}_n\setminus\{h_0\}}\Big|\frac{\widehat{T}_h-\widehat{T}_{h_0}}{\widehat{v}_{h,h_0}}\Big|\ge\gamma_n\Big) \\
& \le\mathbb{P}_n\Big(\max_{h\in\mathcal{H}_n\setminus\{h_0\}}\Big|\frac{\varepsilon'(W_h-W_{h_0})\varepsilon}{v_{h,h_0}}\Big|\ge\frac{\gamma_n}{1+\eta/2}\Big)+o_{\mathbb{P}}(1).
\end{aligned}
$$

Applying Lemma 2(ii) using Assumption W0(iii) and $h_j=h_0a^{-j}$ for $j=0,\dots,J_n$, we obtain

$$
\mathbb{P}_n(\widetilde{h}\neq h_0)\le\sum_{h\in\mathcal{H}_n\setminus\{h_0\}}\mathbb{P}_n\Big(\Big|\frac{\varepsilon'(W_h-W_{h_0})\varepsilon}{v_{h,h_0}}\Big|\ge\frac{\gamma_n}{1+\eta/2}\Big)+o_{\mathbb{P}}(1)
$$

$$\leq \frac{\sqrt{2}(1+\eta/2)}{\sqrt{\pi}\gamma_n} \exp\left(-\frac{1}{2}\left(\frac{\gamma_n}{1+\eta/2}\right)^2 + \ln J_n\right)$$

$$+ O_{\mathbb{P}}(h_0^{p/8}) \sum_{j=1}^{+\infty} \frac{1}{(a^{pj}-1)^{1/8}} + o_{\mathbb{P}}(1) = o_{\mathbb{P}}(1),$$

using (3.8), $h_0 \to 0$ and $\gamma_n \to \infty$. Thus, $\mathbb{P}_n(\widehat{T}_{\widehat{h}} \geq \widehat{v}_{h_0} z_\alpha) = \mathbb{P}_n(\widehat{T}_{h_0} \geq \widehat{v}_{h_0} z_\alpha) + o_{\mathbb{P}}(1)$. Theorem 1 then follows from Propositions 3(ii) and 4, Lemma 2(i) and Assumption W0.

6.5. *Proof of Theorem 4 under Assumptions* D *and* W0. Let $\varepsilon^* = [\varepsilon_1^*, \ldots, \varepsilon_n^*]$. We first establish a moment bound that plays the role of Assumption E. As $\varepsilon_i^* = \widehat{\sigma}_n(X_i)\omega_i$, where the $\omega_i$ are independent of the initial sample, $\mathbb{E}[|\varepsilon_i^*|^{d'}|X_1, Y_1, \ldots, X_n, Y_n] = \mathbb{E}[|\omega_1|^{d'}]|\widehat{\sigma}_n(X_i)|^{d'}$ and

$$(6.8) \quad \max_{1 \leq i \leq n} \mathbb{E}[|\varepsilon_i^*|^{d'}|X_1, Y_1, \ldots, X_n, Y_n] \leq \mathbb{E}[|\omega_1|^{d'}]\left(\sup_{x \in [0,1]^p} \sigma^{d'}(x) + o_{\mathbb{P}}(1)\right).$$

This is sufficient to establish Theorem 4; see Guerre and Lavergne (2003).

6.6. *Proof of Theorem 2 under Assumptions* W0 *and* W1.

LEMMA 6. *Consider a function $\widehat{\delta}(\cdot) \in C_p(L,s)$ with $s > \underline{s}$ and $L > 0$. Consider any sequence $h_n$ from $\mathcal{H}_n$ and let $\Lambda_n = \Lambda_n(s, h_n)$ be as in Assumption* W1(iii). *Under Assumption* W1, *we have*

$$\widehat{\delta}(X)'W_{h_n}\widehat{\delta}(X)$$
$$\geq n[(\Lambda_n - \mathrm{Sp}_n^{1/2}[W_{h_n} - P_{h_n}])\|\widehat{\delta}(X_i)\|_n - (\Lambda_n + \mathrm{Sp}_n^{1/2}[P_{h_n}])C_4 L h_n^s]^2,$$

*where $C_4 = C_4(s)$ is from Assumption* W1(ii), *provided*

$$(6.9) \quad \|\widehat{\delta}(X_i)\|_n \geq \frac{\Lambda_n + \mathrm{Sp}_n^{1/2}[P_{h_n}]}{\Lambda_n - \mathrm{Sp}_n^{1/2}[W_{h_n} - P_{h_n}]} C_4 L h_n^s \geq 0.$$

PROOF. We have $\widehat{\delta}'W_{h_n}\widehat{\delta} = \widehat{\delta}'P_{h_n}\widehat{\delta} + \widehat{\delta}'(W_{h_n} - P_{h_n})\widehat{\delta} \geq \widehat{\delta}'P_{h_n}\widehat{\delta} - n\mathrm{Sp}_n[W_{h_n} - P_{h_n}]\|\widehat{\delta}\|_n^2$. Let $\pi(\cdot)$ be such that $\sup_{x \in [0,1]^p} |\widehat{\delta}(x) - \pi(x)| \leq C_4 L h_n^s$; see Assumption W1(ii). Because $P_{h_n}$ is positive by Assumption W1(i), the triangle inequality and the definition of $\Lambda_n$ yield

$$\left(\frac{\widehat{\delta}'P_{h_n}\widehat{\delta}}{n}\right)^{1/2} \geq \left(\frac{\pi'P_{h_n}\pi}{n}\right)^{1/2} - \left(\frac{1}{n}(\widehat{\delta}-\pi)'P_{h_n}(\widehat{\delta}-\pi)\right)^{1/2}$$

$$\geq \left(\frac{\pi'P_{h_n}\pi}{n}\right)^{1/2} - \mathrm{Sp}_n^{1/2}[P_{h_n}]\|\widehat{\delta}-\pi\|_n$$

$$\geq \Lambda_n \|\widehat{\delta} + \pi - \widehat{\delta}\|_n - \mathrm{Sp}_n^{1/2}[P_{h_n}]\|\widehat{\delta} - \pi\|_n$$

$$\geq \Lambda_n \|\widehat{\delta}\|_n - (\Lambda_n + \mathrm{Sp}_n^{1/2}[P_{h_n}])\|\widehat{\delta} - \pi\|_n$$

$$\geq \Lambda_n \|\widehat{\delta}\|_n - (\Lambda_n + \mathrm{Sp}_n^{1/2}[P_{h_n}])C_4 L h_n^s.$$

As $(\Lambda_n - \mathrm{Sp}_n^{1/2}[W_{h_n} - P_{h_n}])\|\widehat{\delta}\|_n - (\Lambda_n + \mathrm{Sp}_n^{1/2}[P_{h_n}])C_4 L h_n^s \geq 0$ from (6.9),

$$\frac{\widehat{\delta}' W_{h_n} \widehat{\delta}}{n} \geq [\Lambda_n \|\widehat{\delta}\|_n - (\Lambda_n + \mathrm{Sp}_n^{1/2}[P_{h_n}])C_4 L h_n^s]^2 - \mathrm{Sp}_n[W_{h_n} - P_{h_n}]\|\widehat{\delta}\|_n^2$$

$$= [(\Lambda_n - \mathrm{Sp}_n^{1/2}[W_{h_n} - P_{h_n}])\|\widehat{\delta}\|_n - (\Lambda_n + \mathrm{Sp}_n^{1/2}[P_{h_n}])C_4 L h_n^s]$$

$$\times [(\Lambda_n + \mathrm{Sp}_n^{1/2}[W_{h_n} - P_{h_n}])\|\widehat{\delta}\|_n - (\Lambda_n + \mathrm{Sp}_n^{1/2}[P_{h_n}])C_4 L h_n^s]$$

$$\geq [(\Lambda_n - \mathrm{Sp}_n^{1/2}[W_{h_n} - P_{h_n}])\|\widehat{\delta}\|_n - (\Lambda_n + \mathrm{Sp}_n^{1/2}[P_{h_n}])C_4 L h_n^s]^2. \quad \square$$

We now prove Theorem 2 under Assumptions W0 and W1, using the power bound (2.3). Take $h_n = h_0 a^{-j_n}$, where $j_n$ is the integer part of

$$\frac{1}{\ln a}\left[\frac{2}{4s+p}\ln\left(\frac{L^2 n}{\gamma_n \inf_{x\in[0,1]^p}\sigma^2(x)}\right) + \ln h_0\right]$$

$$\asymp \frac{1}{\ln a}\frac{2}{4s+p}\ln\left(\frac{L^2 n}{\gamma_n \inf_{x\in[0,1]^p}\sigma^2(x)}\right),$$

using $\ln h_0 = O(\ln\ln n)$ and $\ln(n/\gamma_n) \geq (1-\gamma)\ln n$ for some $\gamma \in (0,1)$. Note that $h_n$ is in $\mathcal{H}_n$ for all $s > \underline{s}$ and $L > 0$ since $h_{J_n} \asymp (\ln n)^{C_2/p}/n^{2/(4\underline{s}+p)}$ for some $C_2 > 1$ and $\gamma_n \leq n^\gamma$ for some $\gamma \in (0,1)$. We have

$$L h_n^s \asymp L^{p/(4s+p)}\left(\frac{\overline{\sigma}^2 \gamma_n}{n}\right)^{2s/(4s+p)}$$

and

$$n L^2 h_n^{2s} \asymp \gamma_n \overline{\sigma}^2 h_n^{-p/2} \asymp L^{2p/(4s+p)}(\overline{\sigma}^2 \gamma_n)^{4s/(4s+p)} n^{p/(4s+p)} \to \infty.$$

Take now $\widehat{\delta}(\cdot) = m_n(\cdot) - \mu(\cdot; \widehat{\theta}_n)$ in Lemma 6, which belongs to $C_p(L, s)$ by the assumptions of Theorem 2. The lower bound (3.10) of Theorem 2 yields

$$\|\widehat{\delta}(X)\|_n \geq \|m_n(X) - \mu(X; \theta_n)\|_n \geq C\kappa_1 L h_n^s (1 + o_{\mathbb{P}}(1)),$$

implying, in particular, that $n\|m_n(X) - \mu(X; \theta_n)\|_n^2$ diverges in probability. Under Assumptions W0(ii) and W1(i), (iii),

$$\mathbb{P}\left(C\kappa_1 L h_n^s \geq \frac{\Lambda_n(s, h_n) + \mathrm{Sp}_n^{1/2}[P_{h_n}]}{\Lambda_n(s, h_n) - \mathrm{Sp}_n^{1/2}[W_{h_n} - P_{h_n}]}C_4 L h_n^s \geq 0\right) \to 1$$

for $\kappa_1$ large enough, showing that $\widehat{\delta}(\cdot)$ verifies (6.9) with probability tending to 1. Therefore, Lemma 6 and Assumption W1(iii) yield

$$
\begin{aligned}
(m_n(X) &- \mu(X; \widehat{\theta}_n))' W_{h_n}(m_n(X) - \mu(X; \widehat{\theta}_n)) \\
&= \widehat{\delta}'(X) W_{h_n} \widehat{\delta}(X) \\
&\geq n[(\Lambda_n - \mathrm{Sp}_n^{1/2}[W_{h_n} - P_{h_n}]) \| m_n(X) - \mu(X; \theta_n) \|_n \\
&\qquad\qquad - (\Lambda_n + \mathrm{Sp}_n^{1/2}[P_{h_n}]) C_4 L h_n^s]^2 (1 + o_{\mathbb{P}}(1)) \\
&\geq C(1 + o_{\mathbb{P}}(1)) n \| m_n(X) - \mu(X; \theta_n) \|_n^2 \geq C(1 + o_{\mathbb{P}}(1)) n \kappa_1^2 L^2 h_n^{2s}.
\end{aligned}
$$

Moreover, by Proposition 4,

$$
\begin{aligned}
(m_n(X) - \mu(X; \widehat{\theta}_n))' W_{h_n} \varepsilon &= O_{\mathbb{P}}(\sqrt{n} \| m_n(X) - \mu(X; \theta_n) \|_n) \\
&= o_{\mathbb{P}}(n \| m_n(X) - \mu(X; \theta_n) \|_n^2).
\end{aligned}
$$

From $\varepsilon' W_{h_n} \varepsilon = O_{\mathbb{P}}(v_{h_n}) = O_{\mathbb{P}}(h_n^{-p/2}) = o_{\mathbb{P}}(n L^2 h_n^{2s})$ and (6.5),

$$
\widehat{T}_{h_n} \geq C(1 + o_{\mathbb{P}}(1)) n \| m_n(X) - \mu(X; \theta_n) \|_n^2 \geq C(1 + o_{\mathbb{P}}(1)) n \kappa_1^2 L^2 h_n^{2s}.
$$

Proposition 3(ii), Lemma 1 and Assumption W0(iii) yield $z_\alpha \widehat{v}_{h_0} + \gamma_n \widehat{v}_{h_n, h_0} \asymp_{\mathbb{P}} \gamma_n \widehat{v}_{h_n, h_0} \asymp_{\mathbb{P}} \gamma_n \overline{\sigma}^2 h_n^{-p/2} \asymp n L^2 h_n^{2s}$. Collecting the leading terms implies that, for $\kappa_1$ large enough,

$$
\widehat{T}_{h_n} - z_\alpha \widehat{v}_{h_0} - \gamma_n \widehat{v}_{h_n, h_0} \geq C n L^2 h_n^{2s} (\kappa_1^2 - C')(1 + o_{\mathbb{P}}(1)) \xrightarrow{\mathbb{P}} +\infty.
$$

6.7. *Proof of Theorem* 3 *under Assumptions* W0 *and* W1. The proof follows the lines of the proof of Theorem 2, using now (2.4). Since $m_n(X) - \mu(X; \widehat{\theta}_n) = r_n \delta_n(X) + \mu(X; \theta_0) - \mu(X; \widehat{\theta}_n)$,

$$
\begin{aligned}
(m_n(X) &- \mu(X; \widehat{\theta}_n))' W_{h_0}(m_n(X) - \mu(X; \widehat{\theta}_n)) \\
&= r_n^2 \delta_n(X)' W_{h_0} \delta_n(X) \\
&\quad + 2 r_n \delta_n(X) W_{h_0}(\mu(X; \theta_0) - \mu(X; \widehat{\theta}_n)) \\
&\quad + (\mu(X; \theta_0) - \mu(X; \widehat{\theta}_n))' W_{h_0}(\mu(X; \theta_0) - \mu(X; \widehat{\theta}_n)).
\end{aligned}
$$

By Lemma 5,

$$
\begin{aligned}
|r_n \delta_n(X) &W_{h_0}(\mu(X; \theta_0) - \mu(X; \widehat{\theta}_n))| \\
&\leq n r_n \mathrm{Sp}_n[W_{h_0}] \| \delta_n(X) \|_n \| \mu(X; \theta_0) - \mu(X; \widehat{\theta}_n) \|_n = o_{\mathbb{P}}(n r_n^2), \\
|(\mu(X; \theta_0) &- \mu(X; \widehat{\theta}_n))' W_{h_0}(\mu(X; \theta_0) - \mu(X; \widehat{\theta}_n))| \\
&\leq n \mathrm{Sp}_n[W_{h_0}] \| \mu(X; \theta_0) - \mu(X; \widehat{\theta}_n) \|_n^2 = o_{\mathbb{P}}(n r_n^2).
\end{aligned}
$$

Because $\{\delta_n(\cdot)\}_{n\geq 1} \subset C(L, s)$ with $s > \underline{s}$, Lemma 6 yields, under (3.11) and $h_0 \to 0$,

$$
\begin{aligned}
\delta_n(X)'W_{h_0}\delta_n(X) &\geq (1+o_{\mathbb{P}}(1))n[(\Lambda_n - \mathrm{Sp}_n^{1/2}[W_{h_0} - P_{h_0}])\|\delta_n(X)\|_n \\
&\qquad\qquad - C_4(\Lambda_n + \mathrm{Sp}_n^{1/2}[P_{h_0}])Lh_0^s]^2 \\
&\geq Cn(1+o_{\mathbb{P}}(1)).
\end{aligned}
$$

Equation (6.5) in the proof of Proposition 4 and Lemma 5 give, since $z_\alpha\widehat{v}_{h_0} + \varepsilon'W_{h_0}\varepsilon = O_{\mathbb{P}}(h_0^{-p/2})$, $nr_n^2 h_0^{p/2} \to +\infty$ and $h_0 \to 0$,

$$
\widehat{T}_{h_0} - z_\alpha\widehat{v}_{h_0} - \gamma_n\widehat{v}_{h_0,h_0} \geq (1+o_{\mathbb{P}}(1))Cnr_n^2 + O_{\mathbb{P}}(h_0^{-p/2}) \xrightarrow{\mathbb{P}} +\infty.
$$

6.8. *Proof of Proposition* 1. We only detail the case of Examples 1a and 1b. The proof of Proposition 1 for Example 2 can be found in Guerre and Lavergne (2003).

The functions $\psi_k(\cdot)$ can be changed into any system generating the same linear subspace of $\mathbb{R}^n$: Consider the following orthonormal basis of $L_2([0,1]^p, dx)$:

$$
\begin{aligned}
\phi_k(x) &= \prod_{\ell=1}^{p}\sqrt{2k_\ell + 1}\, Q_{k_\ell}(x_\ell)\mathbb{I}(x \in [0,1]^p) \qquad \text{for Example 1a,} \\
(6.10) \quad \phi_{qkh}(x) &= h^{-p/2}\prod_{\ell=1}^{p}\sqrt{2k_\ell + 1}\, Q_{q_\ell}(k_\ell h - x_\ell)\mathbb{I}(x \in I_k(h)) \\
&\qquad\qquad\qquad\qquad\qquad\qquad\qquad \text{for Example 1b,}
\end{aligned}
$$

where the $Q_k(\cdot)$ are the Legendre polynomials of degree $k$ on $[0,1]$, with $\sup_{t\in[0,1]}|Q_k(t)| \leq 1$, $\int_0^1 Q_k^2(t)\,dt = 1/(2k+1)$, $\int_0^1 Q_k(t)Q_{k'}(t)\,dt = 0$ for $k \neq k'$; see, for example, Davis (1975). Let $\Phi_h = [\phi_k(X), 1 \leq |k| \leq 1/h]$ for Example 1a and $\Phi_h = [\phi_{qkh}(X), 1 \leq |q| \leq \bar{q}, 1 \leq |k| \leq 1/h]$ for Example 1b. Define $d_h$ as the number of columns of $\Phi_h$ and note that in both examples $d_h$ is of order $h^{-p}$.

LEMMA 7. *If $f(\cdot)$ is bounded away from 0 and infinity on $[0,1]^p$, there is a $C > 0$ such that*

$$
\max_{h\in\mathcal{H}_n}\mathrm{Sp}_{d_h}[(n^{-1}\Phi_h'\Phi_h)^{-1}] \leq C
$$

*and*

$$
\max_{h\in\mathcal{H}_n}\mathrm{Sp}_{d_h}[n^{-1}\Phi_h'\Phi_h] \leq C \qquad \text{with probability tending to 1,}
$$

*provided $h_{J_n}^{-p} = o(n/\ln n)^{1/3}$ in Example 1a and $h_{J_n}^{-p} = o(n/\ln n)$ in Example 1b.*

PROOF.    Consider first Example 1a. As the $n^{-1}\Phi'_h\Phi_h$, $h \in \mathcal{H}_n$, are nested Gram matrices, it is sufficient to consider the spectral radii of $n^{-1}\Phi'_{h_{J_n}}\Phi_{h_{J_n}}$ and its inverse. We have

$$|\phi_k(X_i)\phi_{k'}(X_i)| \leq \prod_{\ell=1}^{p}\sqrt{2k_\ell+1}\sqrt{2k'_\ell+1} \leq Ch_{J_n}^{-p},$$

$$\mathrm{Var}(\phi_k(X_i)\phi_{k'}(X_i)) \leq \mathbb{E}\phi_k^2(X_i)\phi_{k'}^2(X_i) \leq \mathbb{E}^{1/2}\phi_k^4(X_i)\mathbb{E}^{1/2}\phi_{k'}^4(X_i)$$

$$\leq \sup_{x\in[0,1]^p}|\phi_k(x)|\sup_{x\in[0,1]^p}|\phi_{k'}(x)|\mathbb{E}^{1/2}\phi_k^2(X_i)\mathbb{E}^{1/2}\phi_{k'}^2(X_i)$$

$$\leq Ch_{J_n}^{-p},$$

as $\mathbb{E}\phi_k^2(X) \leq \sup_{x\in[0,1]^p} f(x)\int\phi_k^2(x)\,dx = \sup_{x\in[0,1]^p}f(x)$. The Bernstein inequality then yields

$$\sqrt{\frac{nh_{J_n}^p}{\ln n}}\sup_{0\leq|k|,|k'|\leq 1/h_{J_n}}\left|\frac{1}{n}\sum_{i=1}^n\phi_k(X_i)\phi_{k'}(X_i)-\mathbb{E}\phi_k(X)\phi_{k'}(X)\right| = O_{\mathbb{P}}(1).$$

This gives $n^{-1}\Phi'_{h_{J_n}}\Phi_{h_{J_n}} = n^{-1}\mathbb{E}\Phi'_{h_{J_n}}\Phi_{h_{J_n}} + R_{h_{J_n}}$, where $R_{h_{J_n}}$ is a $d_{h_{J_n}}\times d_{h_{J_n}}$ matrix whose elements are uniformly $O_{\mathbb{P}}(\sqrt{\ln n/nh_{J_n}^p})$. Thus,

$$\mathrm{Sp}_{d_{h_{J_n}}}[R_{h_{J_n}}] \leq N_{d_{h_{J_n}}}[R_{h_{J_n}}] = O_{\mathbb{P}}\left(\frac{1}{h_{J_n}^p}\sqrt{\frac{\ln n}{nh_{J_n}^p}}\right) = o_{\mathbb{P}}(1),$$

as $h_{J_n}^{-p} = o(n/\ln n)^{1/3}$. Hence, the eigenvalues of $n^{-1}\Phi'_{h_{J_n}}\Phi_{h_{J_n}}$ are between the smallest and largest eigenvalues of $n^{-1}\mathbb{E}\Phi'_{h_{J_n}}\Phi_{h_{J_n}}$, with probability tending to one. But, for any $a \in \mathbb{R}^{d_{h_{J_n}}}$,

$$n^{-1}a'\mathbb{E}\Phi'_{h_{J_n}}\Phi_{h_{J_n}}a = \mathbb{E}\left(\sum_{0\leq|k|\leq 1/h_{J_n}}a_k\phi_k(X)\right)^2$$

$$\asymp \int_{[0,1]^p}\left(\sum_{0\leq|k|\leq 1/h_{J_n}}a_k\phi_k(x)\right)^2dx = a'a,$$

since the $\phi_k(\cdot)$ are orthonormal in $L_2([0,1]^p,dx)$. Therefore, the eigenvalues of the symmetric matrix $n^{-1}\mathbb{E}\Phi'_{h_{J_n}}\Phi_{h_{J_n}}$ are bounded away from 0 and infinity when $n$ grows. Example 1b is studied in Baraud (2002) and follows from similar arguments.   □

We now return to the proof of Proposition 1 for Example 1. Lemma 7 implies that, for some $C > 1$,

$$\frac{1}{Cn}\Phi_h\Phi'_h \prec P_h = \frac{1}{n}\Phi_h\left(\frac{1}{n}\Phi'_h\Phi_h\right)^{-1}\Phi'_h \prec \frac{C}{n}\Phi_h\Phi'_h,$$

with probability tending to 1, where $\prec$ is the ordering of symmetric matrices. Because $p_{ii}(h) = e'_i P_h e_i$, where $\{e_i\}_{1 \leq i \leq n}$ is the canonical basis of $\mathbb{R}^n$, this gives

$$(6.11) \quad |p_{ii}(h)| \leq \begin{cases} \dfrac{C}{n} \displaystyle\sum_{|k| \leq 1/h} \phi_k^2(X_i) \leq C/(nh^{2p}), & \text{for Example 1a,} \\[2ex] \dfrac{C}{n} \displaystyle\sum_{|k| \leq 1/h, q \leq \bar{q}} \phi_{qkh}^2(X_i) \leq C/(nh^p), & \text{for Example 1b,} \end{cases}$$

with probability going to 1 and uniformly in $i = 1, \ldots, n$ and $h \in \mathcal{H}_n$. Indeed, $\phi_k^2(\cdot) \leq Ch^{-p}$ for all $k \leq 1/h$ for Example 1a, while $\phi_{qkh}^2(X_i)$ vanishes except for exactly one index $k$ with $\phi_{qkh}^2(X_i) \leq Ch^{-p}$ for Example 1b.

To prove Assumption W0(ii), note that $\mathrm{Sp}_n[P_h] = 1$ since $P_h$ is an orthogonal projection. The triangular inequality gives $\max_{h \in \mathcal{H}_n} \mathrm{Sp}_n[W_h] \leq 1 + \max_{h \in \mathcal{H}_n} \max_{1 \leq i \leq n} |p_{ii}(h)| = O_{\mathbb{P}}(1)$ by (6.11) and the restriction on $h_{J_n}$ which gives $h_{J_n}^{-2p} = o(n)$ for Example 1a and $h_{J_n}^{-p} = o(n)$ for Example 1b. For Assumption W0(iii), we have

$$N_n^2[W_h] = N_n^2[P_h] - N_n^2[W_h - P_h],$$

$$N_n^2[W_h - W_{h_0}] = N_n^2[P_h - P_{h_0}] - N_n^2[(W_h - P_h) - (W_{h_0} - P_{h_0})].$$

Now $N_n^2[P_h] = \mathrm{Rank}[P_h]$ and $N_n^2[P_h - P_{h_0}] = \mathrm{Rank}[P_h - P_{h_0}]$, since $P_h$ and $P_h - P_{h_0}$ are orthogonal projections. This gives $N_n^2[P_h] \asymp h^{-p}$ and $N_n^2[P_h - P_{h_0}] \asymp_{\mathbb{P}} h^{-p} - h_0^{-p}$ almost surely for Example 1a, and for Example 1b, using the Bernstein inequality with $h_{J_n}^{-p} = o(n/\ln n)$, ensuring that the number of $X_i$ in each bin $I_k(h)$ diverges. Then, since $N_n^2[W_h - P_h] = \sum_{i=1}^n p_{ii}^2(h)$, Assumption W0(iii) holds if

$$\max_{h \in \mathcal{H}_n} h^p \sum_{i=1}^n p_{ii}^2(h) = o_{\mathbb{P}}(1)$$

and

$$\max_{h \in \mathcal{H}_n \setminus h_0} (h^{-p} - h_0^{-p})^{-1} \sum_{i=1}^n (p_{ii}(h) - p_{ii}(h_0))^2 = o_{\mathbb{P}}(1),$$

which is a consequence of (6.11), together with $h_{J_n}^{-3p} = o(n/\ln n)$ for Example 1a and $h_{J_n}^{-p} = o(n/\ln n)$ for Example 1b. To show Assumption W1(i), note that the $P_h$ are symmetric positive semidefinite with $\max_{h \in \mathcal{H}_n} \mathrm{Sp}_n[W_h - P_h] = o_{\mathbb{P}}(1)$, as shown when establishing Assumption W0(ii). For Assumption W1(ii), (iii), consider first Example 1a. Let $\Pi_{s,h}$ be the set of polynomial functions with order $1/h$ which are such that Assumption W1(ii) holds by the multivariate Jackson theorem; see, for example, Lorentz (1966). This choice of $\Pi_{s,h}$ gives $\Lambda_n^2 = 1$ almost surely by definition of the $P_h$ with $h_{J_n}^{-p} = o(n)$ and Assumption D. For Example 1b, the proof of Assumtion W1(ii) uses the same Taylor expansion as in Guerre and Lavergne (2002) to build the $\Pi_{s,h}$. Assumption W1(iii), for any given $\bar{q}$, is a consequence of Assumption W1(iii) for $\bar{q} = 1$. This can be shown using Guerre and Lavergne (2002) and establishing convergence of local empirical moments with repeated applications of the Bernstein inequality.

**Acknowledgments.** We thank an Editor, an Associate Editor and three referees for comments that were helpful to improve our paper.

LABORATOIRE DE STATISTIQUE
  THÉORIQUE ET APPLIQUÉE
BOÎTE 158, UNIVERSITÉ PIERRE ET MARIE CURIE
4 PLACE JUSSIEU
75252 PARIS CEDEX 05
FRANCE
E-MAIL: eguerre@ccr.jussieu.fr

GREMAQ
CNRS UMR 5604
UNIVERSITÉ TOULOUSE 1
MANUFACTURE DES TABACS BAT F
21 ALLÉES DE BRIENNE
31000 TOULOUSE
FRANCE
E-MAIL: Pascal.Lavergne@univ-tlse1.fr